\documentclass[11pt]{article}

\usepackage{amssymb,amsmath,amsthm}
\usepackage{graphicx}
\usepackage{multicol}
\usepackage{multirow}

\usepackage{color}

\def\N{\mathbb{N}}
\def\R{\mathbb{R}}
\def\Z{\mathbb{Z}}
\def\i{\mathrm i}
\def\d{\mathrm d}
\def\e{\mathrm e}
\def\E{\mathrm E}
\def\P{\mathrm P}
\def\text{\mbox}

\def\vep{\varepsilon}
\def\1{{\bf 1}}

\newcommand{\noi}{\noindent}
\newcommand {\nn}{\nonumber}

\newtheorem{theorem}{Theorem}[section]
\newtheorem{corollary}
{Corollary}[section]
\newtheorem{lemma}
{Lemma}[section]
\newtheorem{proposition}
{Proposition}[section]

\theoremstyle{definition}

\newtheorem{remark}{Remark}[section]

\definecolor{mdb}{rgb}{0,0.08,0.45}

{\catcode `\@=11 \global\let\AddToReset=\@addtoreset}
\AddToReset{equation}{section}

\begin{document}

\title{Sample covariances of random-coefficient AR(1) panel model}
\author{Remigijus Leipus$^{1}$, \ Anne Philippe$^{2}$, \ Vytaut\.e Pilipauskait\.e$^{3}$, \ Donatas Surgailis$^1$ }
\date{\today \\ \small
\vskip.2cm
$^1$Vilnius University, Faculty of Mathematics and Informatics, Naugarduko 24, 03225 Vilnius, Lithuania\\
$^2$Universit\'{e} de Nantes, Laboratoire de Math\'{e}matiques Jean Leray, 44322 Nantes Cedex 3, France\\
$^3$Aarhus University, Department of Mathematics, Ny Munkegade 118, 8000 Aarhus C, Denmark}

\maketitle

\begin{abstract}

The present paper obtains a complete description of the limit distributions of sample
covariances in $N \times n$ panel data when $N$ and $n$ jointly increase, possibly at different rate.
The panel is formed by $N$ independent samples of length $n$ from random-coefficient AR(1) process with the tail distribution function
of the random coefficient regularly varying at the unit root with exponent $\beta >0$.
We show that for $\beta \in (0, 2)$
the sample covariances may display a variety of stable and non-stable limit behaviors with stability parameter
depending on $\beta$ and the mutual increase rate of $N$ and $n $.

\end{abstract}


\medskip


\section{Introduction}

Dynamic panels providing information on a large population of heterogeneous individuals such as households, firms, etc.
observed at regular time periods, are often described by simple autoregressive models with random parameters near
unity. One of the simplest models for individual evolution is
the random-coefficient AR(1) (RCAR(1)) process
\begin{equation}\label{eq:ar1}
X(t) = a X(t-1) + \vep (t), \quad t \in \Z,
\end{equation}
with standardized i.i.d. innovations $\{\vep(t), t \in \Z \}$ and
a random autoregressive coefficient $a \in [0,1)$
independent of $\{\vep (t), t \in \Z\}$. Granger \cite{gran1980} observed that
in the case when the distribution of $a$ is sufficiently dense near unity
the stationary solution of RCAR(1) equation
in \eqref{eq:ar1} may have
long memory in the sense that the sum of its lagged covariances diverges.
To be more specific, assume that the random coefficient
$a \in [0,1)$ has a density function of the following form
\begin{equation}\label{a:beta}
\phi (x) = \psi(x) (1-x)^{\beta - 1}, \quad x \in [0,1),
\end{equation}
where $\beta >0$ and $\psi(x)$, $x \in [0,1)$ is a bounded function with $\lim_{x \uparrow 1} \psi (x) =: \psi (1) > 0$.
Then for $\beta > 1 $
the covariance function of
stationary solution of RCAR(1) equation
in \eqref{eq:ar1} with
standardized finite variance innovations decays as $ t^{-(\beta -1)} $, viz.,
\begin{equation}\label{gammat}
\gamma (t) := \E X (0) X(t) = \E \frac{a^{|t|}}{1-a^2} \sim (\psi(1)/2) \Gamma(\beta-1) t^{-(\beta-1)}, \qquad
t \to \infty,
\end{equation}
implying
$\sum_{t\in \Z} |\operatorname{Cov} (X(0), X(t))| = \infty$ for $\beta \in (1,2]$.
The same long memory property applies to the contemporaneous aggregate
of $N$ independent individual evolutions $\{X_i(t) \}, i=1, \dots, N$ of \eqref{eq:ar1} and the limit Gaussian
aggregated process arising when $N \to \infty $.
Various properties of the RCAR(1) and more
general RCAR equations were studied in
Gon{\c c}alves and Gouri{\'e}roux \cite{gonc1988}, Zaffaroni \cite{zaff2004},
Celov et al.\ \cite{celo2007},
Oppenheim and Viano \cite{OV2004}, Puplinskait\.e and Surgailis \cite{ps2010},
Philippe et al. \cite{PPS2014} and other works, see Leipus et al.\ \cite{lei2014} for review.

Statistical inference in the RCAR(1) model was discussed in several works.
Leipus et al. \cite{leip2006}, Celov et al. \cite{celo2010} discussed nonparametric estimation of the mixing density $\phi(x)$
using empirical covariances of the limit aggregated process.
For panel RCAR(1) data, Robinson~\cite{bib:ROB78} and
Beran et al.~\cite{ber2010} discussed parametric estimation of the mixing density.
In nonparametric context,
Leipus et al.\ \cite{lei2016} studied estimation of the empirical d.f. of $a$ from panel RCAR(1) observations
and derived its asymptotic properties as $N, n \to \infty$,
while \cite{lei2019} discussed estimation of
$\beta$ in \eqref{a:beta} and testing for long memory in the above panel model.
For a $N \times n$ panel comprising $N$ samples $\{X_i(t), t=1, \dots, n \}$ of length $n$, $i=1, \dots, N$
of independent RCAR(1) processes in \eqref{eq:ar1} with mixing distribution in \eqref{a:beta},
Pilipauskait\.e and Surgailis
\cite{pils2014} studied the asymptotic distribution of the sample mean
\begin{eqnarray}\label{Xmean}
\bar X_{N,n}&:=&\frac{1}{Nn} \sum_{i=1}^{N} \sum_{t=1}^{n} X_i(t)
\end{eqnarray}
as $N, n \to \infty$, possibly at a different rate.
\cite{pils2014}
showed that for $0< \beta < 2 $
the limit distribution of this statistic
depends on whether $N /n^\beta \to \infty$ or
$N /n^\beta \to 0$ in which cases $\bar X_{N,n}$ is asymptotically stable with stability parameter depending on $\beta $ and
taking values in the interval $(0,2]$.
See Table 2 below. As shown in \cite{pils2014}, under the
`intermediate' scaling $N /n^\beta \to c \in (0,\infty)$ the limit distribution of $\bar X_{N,n}$ is more complicated and is given by
a stochastic integral with respect to a certain Poisson random measure.

The present paper discusses asymptotic distribution of sample covariances (covariance estimates)
\begin{eqnarray}\label{covts}
\widehat \gamma_{N,n}(t,s)&:=&\frac{1}{Nn} \sum_{1 \le i, i+s \le N} \sum_{1 \le k, k+t \le n} (X_i(k) - \bar X_{N,n})
(X_{i+s}(k+t)- \bar X_{N,n}), \qquad (t,s) \in \Z^2,
\end{eqnarray}
computed from a similar RCAR(1) panel $\{ X_i (t), \, t = 1, \dots, n, \, i = 1, \dots, N \}$ as in \cite{pils2014},
as $N, n $ jointly increase, possibly at a different rate, and
the lag $(t,s) \in \Z^2 $ is fixed, albeit arbitrary.
Particularly, for $(t,s) = (0,0)$, \eqref{covts} agrees with the sample variance:
\begin{eqnarray}\label{cov00}
{\widehat \gamma_{N,n} (0,0)} = \frac{1}{Nn} \sum_{i=1}^{N} \sum_{k=1}^{n} (X_i(k) - \bar X_{N,n})^2.
\end{eqnarray}
The true covariance function $\gamma(t,s) := \E X_i(k) X_{i+s}(k+t) $ of the RCAR(1) panel model with mixing density in \eqref{a:beta}
exists when $\beta > 1 $ and is given by
\begin{eqnarray} \label{gammats}
&\gamma(t,s) = \begin{cases} \gamma(t), &s=0, \\
0, &s\ne 0,
\end{cases}
\end{eqnarray}
where $\gamma(t)$ defined in \eqref{gammat}. Note that $\gamma(t) $ cannot be recovered from
a single realization of the nonergodic RCAR(1) process $\{X(t)\}$ in \eqref{eq:ar1}. However, the covariance function
in \eqref{gammats} can be consistently estimated from the RCAR(1) $N \times n$ panel when
$N, n \to \infty $, together with rates. The limit distribution of the sample covariance
may exist even for $0< \beta < 1 $ when the covariance itself is undefined.
As it turns out, the limit distribution of $\widehat \gamma_{N,n}(t,s)$ depends on the mutual increase rate of
$N$ and $n$, and is also different for temporal, or iso-sectional lags  ($s = 0$) and cross-sectional lags ($s \ne 0$). The distinctions between the cases $s = 0$ and $s\ne 0$ are due to the fact that, in the latter case,
the statistic in \eqref{covts} involves products $X_i(k)X_{i+s}(k+t)$ of independent processes $X_i $ and $X_{i+s} $,
whereas in the former case, $X_i(k)$ and $X_{i}(k+t)$ are dependent r.v.s.
The main results
of this paper are summarized in Table 1 below. Rigorous formulations are given in Sections 3 and 4.
For better comparison, Table 2 presents the results of \cite{pils2014} about the sample mean
in \eqref{Xmean} for the same panel model.

\begin{table}[htbp]

\begin{center}
\begin{tabular}{| l | l | l |}
\hline
{Mutual increase rate of $N, n$} & {Parameter region} & Limit distribution\\
\hline
{$N/n^{\beta} \to \infty $} & $0 < \beta < 2, \beta \ne 1$ &asymmetric $\beta$-stable\\
\hline
$N/n^{\beta} \to 0 $ & $0 < \beta < 2, \beta \ne 1$ &asymmetric $\beta$-stable \\
\hline
$N/n^{\beta} \to 
c \in (0,\infty) $ & $0 < \beta < 2, \beta \ne 1$ &`intermediate Poisson' \\
\hline
Arbitrary & $ \beta > 2$ &Gaussian \\
\hline
\end{tabular} \\
\vskip.4cm
a) temporal lags ($s = 0$) \\
\end{center}

\begin{center}
\begin{tabular}{| l | l | l |}
\hline
{Mutual increase rate of $N, n$} & {Parameter region} & Limit distribution\\
\hline
{$N/n^{2\beta} \to \infty $} & $1 < \beta < 3/2$ & Gaussian
\\
\cline{2-3}
& $1/2 <\beta < 1$ & symmetric $(2\beta)$-stable \\
\cline{2-3}
\hline
$N/n^{2\beta} \to 0 $ & $3/4 < \beta < 3/2$ & symmetric $(4\beta/3)$-stable \\
\hline
$N/n^{2\beta} \to 
c \in (0,\infty)$ & $3/4 < \beta < 3/2$ & `intermediate Poisson' \\
\hline
Arbitrary & $\beta > 3/2 $ &Gaussian \\
\hline
\end{tabular} \\
\vskip.4cm
b) cross-sectional
lags ($s \ne 0$) \\
\caption{Limit distribution of sample covariances $\widehat \gamma_{N,n}(t,s)$ in \eqref{covts}
 }\label{tab:1}
\end{center}
\end{table}

\begin{table}[htbp]

\begin{center}
\begin{tabular}{| l | l | l |}
\hline
{Mutual increase rate of $N, n$} & {Parameter region} & Limit distribution\\
\hline
{$N/n^{\beta} \to \infty $} & $1 < \beta < 2$ & Gaussian \\
\cline{2-3}
& $0 <\beta < 1$ & symmetric $(2\beta)$-stable \\
\hline
$N/n^{\beta} \to 0 $ & $0 < \beta < 2$ & symmetric $\beta$-stable \\
\hline
$N/n^{\beta} \to 
c \in (0,\infty)$ & $0 < \beta < 2$ & `intermediate Poisson' \\
\hline
Arbitrary  & $\beta > 2 $ & Gaussian \\
\hline
\end{tabular} \\
\vskip.4cm
\caption{Limit distribution of the sample mean $\bar X_{N,n}$ in \eqref{Xmean}
 }\label{tab:2}
\end{center}
\end{table}

\begin{remark}
(i) $\beta$-stable limits in Table 1 a) arising when $N/n^\beta \to 0 $ and $N/n^\beta \to \infty $ have
{\it different } scale parameters and hence the limit distribution of temporal sample covariances is different
in the two cases.

\medskip

\noi (ii) `Intermediate Poisson' limits in Tables 1--2 refer to infinitely divisible distributions defined through certain stochastic integrals w.r.t. Poisson random measure. A similar terminology was used in \cite{pils2016}.

\medskip

\noi (iii) It follows from our results (see Theorem \ref{thmjoint2} below)
that a scaling transition similar as in the case of the sample
mean \cite{pils2014} arises in the interval
$0< \beta < 2 $ for temporal sample covariances and product random fields $X_v (u) X_v(u+t), (u,v) \in \Z^2 $ involving
temporal lags, with the critical rate $N \sim n^\beta $
separating regimes with different limit distributions. For
`cross-sectional' product fields $X_v (u) X_{v+s}(u+t), (u,v) \in \Z^2, s\ne 0  $ involving cross-sectional lags,
a similar scaling transition occurs in the interval $0 < \beta < 3/2 $ with the critical
rate $N \sim n^{2\beta}$ between different scaling regimes, see Theorem \ref{thmjoint}.
The notion of scaling transition for long-range dependent
random fields in $\Z^2$ was discussed in Puplinskait\.e and Surgailis \cite{ps2015}, \cite{ps2016}, Pilipauskait\.e and Surgailis
\cite{pils2016}, \cite{pils2017}.

\medskip

\noi (iv) The limit distributions of cross-sectional sample covariances
in the missing intervals $0< \beta < 1/2 $ and
$0< \beta < 3/4 $ of Table 1 b) are given in Corollary \ref{cor1} below. They are more complicated and
not included in Table 1 b) since the term
$ N n (\bar X_{N,n})^2 $ due to the centering by the sample mean in \eqref{covts} may play the dominating role.

\medskip

\noi (v) We expect that the asymptotic distribution of sample covariances
in the RCAR(1) panel model with common innovations (see \cite{pils2015})
can be analyzed in a similar fashion. Due to the differences between the two
models (the common and the idiosyncratic innovation cases), the asymptotic behavior
of sample covariances might be quite different in these two cases.

\medskip

\noi (vi) The results in Table 1 a) are obtained under the finite 4th moment conditions on the innovations, see
Theorems  \ref{thmjoint2} and \ref{thmclt2} below.
Although the last condition does not guarantee the existence of the 4th moment of the RCAR(1) process,
it is crucial for the limit results, including the CLT in the case $\beta > 2 $.
Scaling transition for sample variances of long-range dependent Gaussian and linear random fields on $\Z^2$ with
finite 4th moment was established in Pilipauskait\.e and Surgailis \cite{pils2017}. On the other side,
Surgailis \cite{sur2004}, Horv\'{a}th and Kokoszka \cite{hor2008} obtained stable limits of sample variances and autocovariances
for long memory moving averages with finite 2nd moment and infinite 4th moment.
Finally, we mention the important works of Davis and Resnick \cite{davis1986} and Davis and Mikosch \cite{davis-mikosch1998} on
limit theory for sample covariance and correlation functions of
moving averages and some nonlinear processes with infinite variance, respectively.

\end{remark}

The rest of the paper is organized as follows. Section 2 presents some preliminary facts, including the
definition and properties of the intermediate processes appearing in Table 1.
Section 3 contains rigorous formulations and the proofs
of the asymptotic results for cross-sectional sample
covariances \eqref{covts}, $s \ne 0$
and the corresponding partial sums processes. Analogous results
for temporal sample covariances and partial sums processes are presented in Section 4.
{Section 4 also contains some applications of these results
to estimation of the autocovariance function $\gamma (t)$ in  \eqref{gammat}
from panel data. }
Some auxiliary proofs are given in Appendix.

\section{Preliminaries}

This section contains some preliminary facts which will be used in the following sections.

\smallskip

\noi {\bf \large 2.1. Double stochastic integrals and quadratic forms.}
Let $B_i, i=1,2 $ be independent standard Brownian motions (BMs) on the real line. Let
\begin{equation}\label{Iij}
I_i(f) := \int_{\R} f(s) \d B_i(s), \qquad
I_{ij}(g) := \int_{\R^2} g(s_1,s_2) \d B_i(s_1) \d B_j(s_2), \quad i,j=1,2,
\end{equation}
denote It\^o-Wiener stochastic integrals (single and double)
w.r.t.\ $B_i, B_j$. 
The integrals in \eqref{Iij}
are jointly defined
for any (non-random) integrands $f \in L^2 (\R),
g \in L^2(\R^2) $;
moreover, $ \E I_i(f) = \E I_{ij} (g) = 0$ and
\begin{eqnarray} \label{covI}
\E I_i(f) I_{i'}(f') &=&\begin{cases}0, &i \ne i', \\
\langle f, f' \rangle, &i=i',
\end{cases}
\quad f, f' \in L^2(\R), \\
\E I_i(f) I_{i'j'}(g)&=&0, \hskip1.5cm \forall i, i', j', \quad f \in L^2(\R), g \in L^2(\R^2), \nn \\
\E I_{ij}(g) I_{i'j'} (g')
&=&\begin{cases}0, &(i,j) \notin \{ (i',j'), (j',i')\}, \\
\langle g, g' \rangle, &(i,j) \in \{ (i',j'), (j',i')\}, \ i\ne j, \\
2 \langle g, {\rm sym} g'\rangle, &i = i' = j = j',
\end{cases}
\qquad g, g' \in L^2(\R^2), \nn
\end{eqnarray}
where $\langle f, f' \rangle = \int_{\R} f(s) f'(s) {\rm d} s$ $(\| f \| := \sqrt{\langle f, f \rangle}),
\langle g, g' \rangle = \int_{\R^2} g(s_1,s_2) g'(s_1,s_2) \d s_1 \d s_2$ $(\|g \| := \sqrt{\langle g, g \rangle})$ denote scalar products (norms) in
$L^2(\R)$ and $L^2(\R^2)$, respectively, and ${\rm sym}$ denotes the symmetrization,
see, e.g., (\cite{gir2012}, sec.\ 11.5, 14.3). Note that for $g(s_1,s_2) = f_1(s_1) f_2(s_2), f_i \in L^2(\R), i=1,2 $
we have $I_{ii} (g) = I_i(f_1) I_i(f_2) - \langle f_1, f_2 \rangle, \ I_{12} (g) = I_1(f_1) I_2(f_2), $ in particular,
$I_{12} (g) =_{\rm d} \|f_1\| \|f_2\| Z_1 Z_2, $ where $Z_i \sim N(0,1), i=1,2 $ are independent standard normal r.v.s.

Let $\{ \vep_i(s), s \in \Z\}, i=1,2 $ be independent sequences of standardized i.i.d. r.v.s, $\E \vep_i (s) = 0$, $\E \vep_i (s) \vep_{i'} (s') = 1$ if $(i,s) = (i',s')$, $\E \vep_i (s) \vep_{i'} (s') = 0$ if $(i,s) \neq (i',s')$, $i,i' = 1,2$, $s,s' \in \Z$.
Consider the centered quadratic form
\begin{equation}\label{Qij}
Q_{ij} (h) = \sum_{s_1, s_2 \in \Z} h (s_1,s_2) [\vep_i(s_1) \vep_j(s_2) - \E \vep_i(s_1) \vep_j(s_2)], \qquad i,j=1,2,
\end{equation}
where $h \in L^2(\Z^2)$. For $i=j$ we additionally assume $\E \vep^4_i(0) < \infty $.
Then the sum in \eqref{Qij} converges in $L^2 $ and
\begin{equation}
{\rm var} (Q_{ij} (h)) \le (1 + \E \vep^4_i(0) \delta_{ij}) \sum_{s_1, s_2 \in \Z} h^2 (s_1,s_2),
\end{equation}
see (\cite{gir2012}, (4.5.4)). With any $h \in L^2(\Z^2)$ and any $\alpha_1, \alpha_2 >0$ we associate its extension to $L^2(\R^2)$, namely,
\begin{equation}
\widetilde h^{(\alpha_1,\alpha_2)}(s_1,s_2) := (\alpha_1 \alpha_2)^{1/2} h(\lfloor\alpha_1 s_1\rfloor, \lfloor\alpha_2 s_2\rfloor), \qquad (s_1,s_2) \in \R^2,
\end{equation}
with $\| \widetilde h^{(\alpha_1,\alpha_2)}\|^2 = \sum_{s_1,s_2 \in \Z} h^2(s_1,s_2)$.
We shall use the following criterion for the convergence
in distribution of quadratic forms in \eqref{Qij} towards double stochastic integrals
\eqref{Iij}.

\begin{proposition}\label{Qijconv} {\rm (\cite{gir2012}, Proposition\,11.5.5)} Let $i,j=1,2 $ and
$Q_{ij} (h_{\alpha_1,\alpha_2}), \alpha_1, \alpha_2 >0$ be a family of quadratic forms as in \eqref{Qij} with
coefficients $h_{\alpha_1,\alpha_2} \in L^2(\Z^2)$. For $i=j$ we additionally assume $\E \vep^4_i(0) < \infty $.
Suppose for some $g \in L^2(\R^2)$ 
we have that
\begin{equation}
\lim_{\alpha_1, \alpha_2 \to \infty } \| \widetilde h^{(\alpha_1,\alpha_2)}_{\alpha_1,\alpha_2} - g\| = 0.
\end{equation}
Then $Q_{ij} (h_{\alpha_1,\alpha_2}) \to_{\rm d} I_{ij}(g) \ (\alpha_1,\alpha_2 \to \infty)$, where $I_{ij}(g) $ is defined as in  \eqref{Iij}.

\end{proposition}

\medskip

\noi {\bf \large 2.2. The `cross-sectional' intermediate process.} Let $\d {\cal M}_\beta \equiv {\cal M}_\beta(\d x_1, \d x_2, \d B_1, \d B_2)$ denote Poisson random measure
on $(\R_+ \times C(\R))^2$ with mean
\begin{equation} \label{mubeta}
\d \mu_\beta \equiv \mu_\beta (\d x_1, \d x_2, \d B_1, \d B_2)
:= \psi(1)^2 (x_1 x_2)^{\beta -1} \d x_1 \d x_2 P_B(\d B_1) P_B(\d B_2),
\end{equation}
where $\beta >0$ is parameter and
$P_B $ is the Wiener measure on $C(\R)$. Let $ \d\widetilde {\cal M}_\beta := \d {\cal M}_\beta - \d \mu_\beta $ be the centered
Poisson random measure. We shall often use finiteness of the following integrals:
\begin{eqnarray}\label{Aineq2}
&\int_{\R^2_+} \min \big\{1, \frac{1}{x_1 x_2 (x_1 + x_2)}\big\} (x_1 x_2)^{\beta -1} \d x_1 \d x_2 < \infty, \qquad \forall \, 0< \beta < 3/2, \\
&\int_{\R^2_+} \min \big\{1, \frac{1}{x_1 + x_2}\big\} (x_1 x_2)^{\beta -2} \d x_1 \d x_2 \ < \ \infty, \qquad \forall \, 1< \beta < 3/2, \label{Ez5}
\end{eqnarray}
see Appendix. Let
\begin{equation}\label{calY}
{\cal Y}_i(u;x) = \int_{-\infty}^u \e^{-x(u-s) } \d B_i(s), \qquad u \in \R, \ x >0,
\end{equation}
be a family of stationary Ornstein-Uhlenbeck (O-U) processes subordinated to $B_i = \{ B_i (s), s \in \R \}$, $B_i, i=1,2$ being independent BMs.
Let
\begin{eqnarray} \label{zdef}
z(\tau; x_1,x_2)&:=&\int_0^\tau \prod_{i=1}^2 {\cal Y}_i(u;x_i) \d u, \qquad \tau \ge 0, 
\end{eqnarray}
be a family of
integrated products of independent O-U processes indexed by $x_1, x_2 >0$. We use the representation of \eqref{zdef}
\begin{eqnarray} \label{zdefint}
z(\tau; x_1,x_2)&=&\int_{\R^2} \big\{ \int_0^\tau \prod_{i=1}^2 \e^{-x_i(u-s_i)} \1(u>s_i) \d u \big\} \d B_1(s_1) \d B_2(s_2)
\end{eqnarray}
as the double It\^o-Wiener integral in \eqref{Iij}.
The `cross-sectional' intermediate process ${\cal Z}_\beta $ is defined as stochastic integral w.r.t. the Poisson measure
${\cal M}_\beta$, viz.,
\begin{eqnarray}\label{interZ0}
{\cal Z}_\beta(\tau)
&:=&\int_{{\cal L}_1} z(\tau; x_1,x_2) \d {\cal M}_\beta 
+ \int_{{\cal L}^c_1} z(\tau; x_1,x_2) \d \widetilde {\cal M}_\beta, 
\end{eqnarray}
where
\begin{equation}
{\cal L}_1 := \{(x_1,x_2,B_1,B_2)\in (\R_+ \times C(\R))^2: x_1 x_2(x_1+x_2) \le 1 \},
\qquad {\cal L}^c_1 := (\R_+ \times C(\R))^2\setminus {\cal L}_1
\end{equation}
and $\mu_\beta ({\cal L}_1) < \infty $. For $ 1/2 < \beta < 3/2 $ the two integrals
in \eqref{interZ0} can be combined in a single one:
 \begin{eqnarray}\label{interZ}
{\cal Z}_\beta(\tau)
&=&
\int_{(\R_+ \times C(\R))^2} z(\tau; x_1,x_2) \d \widetilde {\cal M}_\beta.
\end{eqnarray}
These and other properties of ${\cal Z}_\beta$ are stated in the following
proposition whose proof is given in the Appendix. We also refer to \cite{rajp1989} and
\cite{pils2014} for general properties of stochastic integrals w.r.t. Poisson random measure.

\begin{proposition} \label{propinter} (i) The process ${\cal Z}_\beta $ in \eqref{interZ0} is well-defined for any $0< \beta < 3/2 $. It has
stationary increments, infinitely divisible finite-dimensional distributions, and the joint ch.f. given by
\begin{eqnarray}\label{Zchf}
\E \exp \big\{ \i\sum_{j=1}^m \theta_j {\cal Z}_\beta (\tau_j) \big\}
&=&\exp \Big\{ \int_{(\R_+ \times C(\R))^2} \big(\e^{\i \sum_{j=1}^m \theta_j z(\tau_j;x_1,x_2)} -1 \big) \d \mu_\beta \Big\},
\end{eqnarray}
where $\theta_j \in \R, \tau_j \ge 0,
j=1, \dots, m, \, m \in \N$. Moreover, the distribution of ${\cal Z}_\beta $ is
symmetric: $\{{\cal Z}_\beta (\tau), \tau \ge 0 \} =_{\rm fdd} \{-{\cal Z}_\beta (\tau), \tau \ge 0 \}.  $

\smallskip

\noi (ii) $\E |{\cal Z}_\beta (\tau)|^p < \infty $ for $p < 2\beta $ and $\E {\cal Z}_\beta (\tau) = 0 $
for $1/2 < \beta < 3/2 $.

\smallskip

\noi (iii) For $1/2 < \beta < 3/2$, ${\cal Z}_\beta $ can be defined as in \eqref{interZ}. Moreover, if $1 < \beta < 3/2$, then
$\E {\cal Z}^2_\beta (\tau) < \infty $ and
\begin{equation} \label{covZ}
\E {\cal Z}_\beta (\tau_1) {\cal Z}_\beta (\tau_2) = (\sigma^2_\infty/2) \big(\tau_1^{2(2-\beta)} + \tau_2^{2(2-\beta)} - |\tau_2 - \tau_1|^{2(2-\beta)}\big),
\quad \tau_1, \tau_2 \ge 0,
\end{equation}
where $\sigma^2_\infty := \psi(1)^2 \Gamma(\beta -1)^2/(4 (2-\beta)(3-2\beta))$.

\smallskip

\noi (iv) For $1/2 < \beta < 3/2 $, the process ${\cal Z}_\beta $ has a.s.\ continuous trajectories.

\smallskip

\noi (v) (Asymptotic self-similarity) As $b \to 0$,
\begin{align}\label{ass:local1}
b^{\beta-2} {\cal Z}_\beta (b \tau) &\to_{\rm fdd} \sigma_\infty B_{2-\beta} (\tau), \qquad \text{if } 1 < \beta < 3/2,\\
b^{-1} (\log b^{-1} )^{-1/(2\beta) } {\cal Z}_\beta (b\tau) &\to_{\rm fdd} \tau V_{2\beta}, \qquad \qquad \ \text{if } 0 < \beta < 1,\label{ass:local0}
\end{align}
where $\{B_{2-\beta}(\tau), \tau \ge 0\} $ is a fractional Brownian motion with $ \E [B_{2-\beta}(\tau)]^2 = \tau^{2(2-\beta)}$, $\tau \ge 0$, $2-\beta \in (1/2, 1)$,
$\sigma^2_\infty$ is given in \eqref{covZ},
and $V_{2\beta}$ is a symmetric $(2\beta)$-stable r.v. with
ch.f. $\E \e^{\i \theta V_{2\beta}} = \e^{-c_\infty |\theta|^{2\beta}}, \theta \in \R$,
$c_\infty := \psi (1)^2 2^{1-2\beta} \Gamma(\beta + (1/2))
\Gamma (1-\beta)/ \sqrt{\pi} $.
For any $0 < \beta < 3/2$, as $b \to \infty$,
\begin{equation}\label{ass:global}
b^{-1/2} {\cal Z}_\beta (b\tau) \to_{\rm fdd} {\cal A}^{1/2} B (\tau),
\end{equation}
where ${\cal A} >0$ is a $(2\beta/3)$-stable r.v.\ with Laplace transform $\E \e^{ - \theta {\cal A} } = \exp \{ - \sigma_0 \theta^{2\beta/3} \}$, $\theta \ge 0$, $\sigma_0 := \psi (1)^2 2^{-2\beta/3} \Gamma(1 - (2\beta/3)) \operatorname{B}(\beta/3, \beta/3)/(2\beta) $,
and $\{ B(\tau), \, \tau \ge 0 \}$ is a standard BM, independent of ${\cal A}$. Finite-dimensional
distributions of the limit process in \eqref{ass:global} are symmetric $(4\beta/3)$-stable.
\end{proposition}

\medskip

\noi {\bf \large 2.3. The `iso-sectional' intermediate process.}
Let $\d {\cal M}^\ast_\beta \equiv {\cal M}^\ast_\beta(\d x, \d B)$ denote Poisson random measure
on $\R_+ \times C(\R)$ with mean
\begin{equation} \label{mubeta1}
\d \mu^\ast_\beta \equiv \mu^\ast_\beta (\d x, \d B)
:= \psi(1) x^{\beta -1} \d x P_B(\d B),
\end{equation}
where $0<\beta <2$ is parameter and
$P_B $ is the Wiener measure on $C(\R)$. Let $ \d\widetilde {\cal M}^\ast_\beta := \d {\cal M}^\ast_\beta - \d \mu^\ast_\beta $ be the centered
Poisson random measure.
Let ${\cal Y}(\cdot;x) \equiv {\cal Y}_1(\cdot;x)$ be the family of O-U processes as in \eqref{calY}, and
\begin{eqnarray} \label{zdef1}
z^\ast (\tau; x)&:=&  
\int_0^\tau {\cal Y}^2(u;x) \d u, \qquad \tau \ge 0, \ x>0,
\end{eqnarray}
be integrated squared O-U processes. Note $\E z^\ast (\tau; x) = \tau \E {\cal Y}^2(0;x) = \tau \int_{-\infty}^0 \e^{2xs} \d s = \tau/(2x)$.
We will use the representation
\begin{equation} \label{zdef1int}
z^\ast (\tau; x)\ =\ \int_{\R^2} \big\{ \int_0^\tau \prod_{i=1}^2 \e^{-x(u-s_i)} \1(u>s_i) \d u \big\} \d B(s_1) \d B(s_2)
+ \tau/(2x)
\end{equation}
as the double It\^o-Wiener integral.
The `iso-sectional' intermediate process ${\cal Z}^\ast_\beta $ is defined for $\beta \in (0,2), \beta \ne 1 $
as stochastic integral w.r.t.\ the above Poisson measure,
viz.,
\begin{eqnarray}\label{interZ01}
{\cal Z}^\ast_\beta(\tau)
&:=&\int_{\R_+ \times C(\R)} z^\ast(\tau; x) \begin{cases}
\d {\cal M}^\ast_\beta, &0< \beta < 1, \\
\d \widetilde {\cal M}^\ast_\beta, &1< \beta < 2,
\end{cases} \qquad \tau \ge 0.
\end{eqnarray}
 Proposition \ref{propinter1} stating properties of ${\cal Z}^\ast_\beta$ is similar to
 Proposition \ref{propinter}.

\begin{proposition} \label{propinter1} (i) The process ${\cal Z}^\ast_\beta $ in \eqref{interZ01} is well-defined for any $0< \beta < 2, \beta \ne 1 $. It has
stationary increments, infinitely divisible finite-dimensional distributions, and the joint ch.f. given by
\begin{eqnarray}\label{Zchf1}
\hskip-1cm &\E \exp \big\{ \i\sum_{j=1}^m \theta_j {\cal Z}^*_\beta (\tau_j) \big\}
=\exp \big\{ \int_{\R_+ \times C(\R)} \big(\e^{\i \sum_{j=1}^m \theta_j z^\ast(\tau_j;x)} -1
- \i \sum_{j=1}^m \theta_j z^\ast(\tau_j;x) \1(1 < \beta < 2)
\big) \d \mu^\ast_\beta \big\},
\end{eqnarray}
where $\theta_j \in \R$, $\tau_j \ge 0$, $j=1, \dots, m$, $m \in \N$.

\smallskip

\noi (ii) $\E |{\cal Z}^\ast_\beta (\tau)|^p < \infty$ for any $0< p < \beta < 2$, $\beta \ne 1 $ and $\E {\cal Z}^\ast_\beta (\tau) = 0$ for $1 < \beta < 2$.



\smallskip

\noi (iii) For $1 < \beta < 2$, the process ${\cal Z}^\ast_\beta$ has a.s.\ continuous trajectories.

\smallskip

\noi (iv) (Asymptotic self-similarity) For any $0 < \beta < 2$, $\beta \neq 1$,
\begin{align}
b^{-1} {\cal Z}_\beta^\ast (b \tau) \to_{\rm fdd}
\begin{cases}
\tau V^\ast_\beta &\text{as } b \to 0,\\
\tau V^+_\beta &\text{as } b \to \infty,
\end{cases}
\end{align}
where $V^+_\beta$, $V^\ast_\beta$ are a completely asymmetric $\beta$-stable r.v.s with ch.f.s $\E \e^{\i \theta V^+_\beta} = \exp \{\psi(1) \int_0^\infty (\e^{\i \theta/(2x)} - 1 - \i (\theta/(2x)) \1(1< \beta < 2) ) x^{\beta -1} \d x \}$,
$\E \e^{\i \theta V^\ast_\beta}
= \exp \{\psi(1) \int_0^\infty \E (\e^{\i \theta Z^2/(2x)}-1-\i (\theta Z^2/(2x)) \1(1< \beta < 2) ) x^{\beta -1} \d x \}$, $\theta \in \R$ and $Z \sim N(0,1)$.
\end{proposition}


\medskip

\noi {\bf \large 2.4. Conditional long-run variance of products of RCAR(1) processes.}
We use some facts in Proposition \ref{propvar}, below,
about conditional variance of the partial sums process
of the product $Y_{ij}(t) := X_i(t) X_j(t) $ of two RCAR(1) processes. Split
$Y_{ij}(t) 
= Y^+_{ij}(t) + Y_{ij}^-(t)$,
where $Y_{ij}^+(t) = \sum_{s_1 \wedge s_2 \ge 1} a_i^{t-s_1} a_j^{t-s_2}\1(t\ge s_1\vee s_2) \vep_i(s_1) \vep_j(s_2),
Y_{ij}^-(t) = \sum_{s_1 \wedge s_2 \le 0} a_i^{t-s_1} a_j^{t-s_2}\1(t\ge s_1\vee s_2) \vep_i(s_1) \vep_j(s_2) $.
For $i=j $ we assume additionally that $ \E \vep_i^4(0) < \infty$.

\medskip

\begin{proposition} \label{propvar}
We have 
\begin{eqnarray} \label{var2}
{\rm var} \big[\sum_{t=1}^n Y_{ij}(t)|a_i, a_j\big]&\sim&{\rm var}\big[\sum_{t=1}^n Y^+_{ij}(t)|a_i, a_j\big] \ \sim \ A_{ij} n, \qquad n \to \infty,
\end{eqnarray}
where
\begin{equation}\label{Aijdef}
A_{ij} := \begin{cases}
\frac{1+ a_i a_j}{(1-a_i^2)(1-a_j^2)(1-a_ia_j)}, &i\ne j, \\
\frac{1+a^2_{i}}{1-a_i^2} (\frac{2}{(1-a_i^2)^2} + \frac{{\rm cum}_4}{1-a_i^4} ), &i=j
\end{cases}
\end{equation}
with ${\rm cum}_4 $ being the 4th cumulant of $\vep_i(0)$.
Moreover, for any  $n\ge 1,\, i,j \in \Z,
\ a_i, a_j \in [0,1)$
\begin{eqnarray}
&&{\rm var} \big[
\sum_{t=1}^n Y_{ij}(t)|a_i,a_j \big] \  \le \  \frac{C_{ij}n^2}{(1-a_i)(1-a_j)}
\min \big\{ 1, \frac{1}{n(2-a_i -a_j)}\big\}, \label{var1}
\end{eqnarray}
where $C_{ij} := 4 \ (i\ne j)$, $:= 2(2+ |{\rm cum}_4|)\ (i=j)$.

\end{proposition}

\noi {\it Proof.} Let $i\ne j$.
We have $\E [Y_{ij}(t) Y_{ij}(s)| a_i, a_j] = \E [X_{i}(t) X_{i}(s)| a_i] \E [X_{j}(t) X_{j}(s)| a_j]
 =  (a_i a_j)^{|t-s|} /(1-a_i^2) (1-a_j^2) $
and hence
\begin{eqnarray} \label{var4}
J_n(a_i,a_j) := \E \big[ \big(\sum_{t=1}^n Y_{ij}(t)\big)^2|a_i,a_j\big]
&=&\frac{n}{(1-a_i^2)(1-a_j^2)}\sum_{t=-n}^n (a_i a_j)^{|t|} \big(1 - \frac{|t|}{n}\big).
\end{eqnarray}
Relation \eqref{var4} implies
\eqref{var2}. It also implies $J_n(a_i,a_j)
\le 2n^2/((1-a_i)(1-a_j))$.
Note also $1 - a_i a_j \ge (1/2)((1-a_i) + (1-a_j))$.
Hence and from \eqref{var4} we obtain
$$
 J_n(a_i,a_j) \le \frac{n}{(1-a_i^2)(1-a_j^2)} \big(1 + 2 \sum_{t=1}^\infty (a_i a_j)^t\big)
\le \frac{2n}{(1-a_i)(1-a_j)(1-a_i a_j)} \le \frac{4n}{(1-a_i)(1-a_j)(2-a_i -a_j)},
$$ proving \eqref{var1}.
The proof of \eqref{var2}--\eqref{var1} for $i=j$ is similar using
${\rm cov}[Y_{ii}(t), Y_{ii}(s)|a_i] = 2 (a_i^{|t-s|}/(1-a_i^2))^2 + {\rm cum}_4 a^{2|t-s|}_i/(1-a^4_i). $
\hfill $\Box$

\medskip

\section{Asymptotic distribution of cross-sectional sample covariances  
}\label{sec3}

Theorems \ref{thmjoint} and \ref{thmclt} discuss the asymptotic distribution of partial sums process
\begin{equation}\label{Sts}
S^{t,s}_{N,n} (\tau) := \sum_{i=1}^{N} \sum_{u=1}^{\lfloor n\tau \rfloor} X_i (u) X_{i+s}(u+t), \quad \tau \ge 0,
\end{equation}
where $t$ and $s \in \Z, s \ne 0 $ are fixed and
$N$ and $n$ tend to infinity, possibly at a different rate.
The asymptotic behavior of sample covariances
$\widehat \gamma_{N,n}(t,s) $ is discussed in Corollary \ref{cor1}.
As it turns out, these limit distributions do not depend
on $t, s $ which is due to the fact that the sectional processes $\{ X_i(t), t \in \Z \}, i \in \Z $ are
independent and stationary.

\begin{theorem} \label{thmjoint} Let the mixing distribution satisfy condition \eqref{a:beta} with $0 < \beta < 3/2 $. Let
$N, n \to \infty $ so as
\begin{equation}\label{mucond}
\lambda_{N,n} := \frac{N^{1/(2\beta)}}{n} \ \to \lambda_\infty \ \in \ [0, \infty].
\end{equation}
Then the following statements (i)--(iii) hold for $S^{t,s}_{N,n}(\tau), (t,s) \in \Z^2, s \ne 0 $ in \eqref{Sts}
depending on $\lambda_\infty$ in \eqref{mucond}.
\medskip

\noi (i) Let $\lambda_\infty = \infty$. Then
\begin{eqnarray} \label{muinf}
&&\hskip1cm 
n^{-2} \lambda_{N,n}^{-\beta}
S^{t,s}_{N,n}(\tau) \quad \to_{\rm fdd} \quad \sigma_\infty B_{2-\beta}(\tau), \qquad  1 < \beta < 3/2, \\
&&
n^{-2} \lambda^{-1}_{N,n} (\log \lambda_{N,n})^{-1/(2\beta)}
S^{t,s}_{N,n}(\tau) \ \to_{\rm fdd}\ \tau V_{2\beta}, \qquad  \ 0< \beta < 1, \label{mulog}
\end{eqnarray}
where the limit processes are the same as in \eqref{ass:local1}, \eqref{ass:local0}.

\medskip

\noi (ii) Let $\lambda_\infty = 0$ and
$\E |\vep(0)|^{2p} < \infty$ for some $p> 1$.
Then
\begin{equation}\label{muzero}
n^{-2} \lambda_{N,n}^{-3/2} S^{t,s}_{N,n}(\tau) \ \to_{\rm fdd} \ {\cal A}^{1/2} B(\tau),
\end{equation}
where the limit process is the same as in \eqref{ass:global}.

\medskip

\noi (iii) Let $0<\lambda_\infty < \infty$. Then
\begin{eqnarray} \label{interlim}
&&
n^{-2} \lambda_{N,n}^{-3/2}
S^{t,s}_{N,n}(\tau)\ \to_{\rm fdd}\ \lambda_\infty^{1/2} {\cal Z}_\beta(\tau/\lambda_\infty),
\end{eqnarray}
where ${\cal Z}_\beta $ is the intermediate process in \eqref{interZ0}.

\end{theorem}

\begin{theorem} \label{thmclt} Let the mixing distribution satisfy condition \eqref{a:beta} with $\beta > 3/2 $ and assume $\E |\vep(0)|^{2p} < \infty$ for some $p>1$. Then for any
$(t,s) \in \Z^2, s\ne 0$
as $N, n \to \infty $ in arbitrary way,
\begin{equation}
n^{-1/2} N^{-1/2} S^{t,s}_{N,n} (\tau) \ \to_{\rm fdd}\ \sigma B(\tau), \qquad \sigma^2 := \E A_{12},
\end{equation}
where $A_{12} $ is defined in \eqref{Aijdef}.

\end{theorem}

\begin{remark}
Our proof of Theorem \ref{thmjoint} (ii) requires establishing the asymptotic normality of a bilinear form in i.i.d. r.v.s, which has a non-zero diagonal, see the r.h.s.\ of \eqref{zzlim}. For this purpose, we use the martingale CLT and impose an additional condition of $\E |\vep (0)|^{2p} < \infty$, $p>1$. To establish the CLT for quadratic forms with non-zero diagonal, \cite{bhan2007} took similar approach and also needed $2p$ finite moments.

In Theorem \ref{thmclt} we also assume $\E |\vep (0)|^{2p} < \infty$, $p>1$. However, it can be proved under $\E \vep^2 (0)< \infty$ applying another technique that is approximation by $m$-dependent r.v.s. Moreover, this result holds if \eqref{a:beta} is replaced by $\E A_{12}<\infty$.

		
\end{remark}

Note that the asymptotic distribution of sample covariances $\widehat \gamma_{N,n} (t,s)$ in \eqref{covts}
coincides with that of the statistics
\begin{equation} \label{tildeG}
\widetilde \gamma_{N,n} (t,s)
:= (Nn)^{-1} S^{t,s}_{N,n} (1) - (\bar X_{N,n})^2.
\end{equation}
For $s\ne 0$ the limit behavior of the first term on the r.h.s.\ of \eqref{tildeG}
can be obtained from Theorems \ref{thmjoint} and \ref{thmclt}. It turns out that for some values
of $\beta $, the second term on the r.h.s.\ can play the dominating role.
The limit behavior of $\bar X_{N,n} $ was identified in \cite{pils2014}
and is given in the following proposition, with some simplifications.

\begin{proposition} \label{propX}
Let the mixing distribution satisfy condition \eqref{a:beta} with $\beta >0$.

\medskip

\noi (i) Let $1< \beta < 2$ and $N/n^\beta \to \infty $. Then
\begin{eqnarray}
N^{1/2} n^{(\beta -1)/2} \bar X_{N,n} \to_{\rm d} \bar \sigma_\beta Z,
\end{eqnarray}
where $Z \sim N(0,1)$ and $\bar \sigma^2_\beta := \psi(1) \Gamma(\beta -1)/((3-\beta)(2-\beta))$.

\medskip

\noi (ii) Let $0< \beta < 1$ and $N/n^\beta \to \infty $. Then
\begin{eqnarray}
N^{1 - 1/2\beta} \bar X_{N,n} \to_{\rm d} \bar V_{2\beta},
\end{eqnarray}
where $\bar V_{2\beta} $ is a symmetric $(2\beta)$-stable r.v. with ch.f. $\E \e^{\i \theta \bar V_{2\beta} }
= \e^{-\bar K_\beta |\theta|^{2\beta}}, \bar K_\beta := \psi(1) 4^{-\beta} \Gamma(1-\beta)/\beta $.

\medskip

\noi (iii) Let $0< \beta < 2$ and $N/n^\beta \to 0 $. Then
\begin{eqnarray}
N^{1- 1/\beta} n^{1/2} \bar X_{N,n} \to_{\rm d} \bar W_\beta,
\end{eqnarray}
where $\bar W_{\beta} $ is a symmetric $\beta$-stable r.v. with ch.f. $\E \e^{\i \theta \bar W_{\beta} }
= \e^{-\bar k_\beta |\theta|^{\beta}}, \bar k_\beta := \psi(1) 2^{-\beta/2} \Gamma(1-\beta/2)/\beta $.

\noi (iv) Let $\beta > 2$. Then as $N, n \to \infty $ in arbitrary way,
\begin{eqnarray}
N^{1/2} n^{1/2} \bar X_{N,n} \to_{\rm d} \bar \sigma Z,
\end{eqnarray}
where $Z \sim N(0,1)$ and $\bar \sigma^2 := \E (1-a)^{-2}$.

\end{proposition}

From Theorems \ref{thmjoint} and Proposition \ref{propX} we see that the r.h.s.\ of
\eqref{tildeG} may exhibit {\it two} `bifurcation points' of the limit behavior, viz.,
as $N \sim n^{2\beta}$ and $N \sim n^\beta$. Depending on the value of $\beta $
the first or the second term may dominate, and the limit behavior of $\widehat \gamma_{N,n}(t,s) $
gets more complicated. The following corollary provides this limit without detailing
the `intermediate' situations and also with exception of some particular values
of $\beta $ where {\it both} terms on the r.h.s.\ may contribute to the limit.
Essentially, the corollary follows by comparing the normalizations
in Theorems \ref{thmjoint} and Proposition \ref{propX}.



\begin{corollary} \label{cor1}
Assume that the mixing distribution satisfies condition \eqref{a:beta} with $\beta >0$ and $\E |\vep(0)|^{2p} < \infty$ for some $p > 1$
and $(t,s) \in \Z^2, s \ne 0$ be fixed albeit arbitrary.

\medskip

\noi (i) Let $N/n^{2\beta} \to \infty $ and $1 < \beta < 3/2 $. Then
\begin{eqnarray*}
&N^{1/2} n^{\beta -1} \widehat \gamma_{N,n} (t,s) \ \to_{\rm d} \ \sigma_\infty Z,
\end{eqnarray*}
where $Z \sim N(0,1)$ and $\sigma_\infty$ is the same as in Theorem \ref{thmjoint} (i).

\medskip

\noi (ii) Let $N/n^{2\beta} \to \infty $ and $1/2 < \beta < 1$. Then
\begin{eqnarray*}
&\frac{N^{1-1/(2\beta)}}{
\log^{1/(2\beta)} (N^{1/(2\beta)}/n)}\, \widehat \gamma_{N,n} (t,s) \to_{\rm d} V_{2\beta},
\end{eqnarray*}
where $V_{2\beta}$ is symmetric $(2\beta)$-stable r.v. defined in Theorem \ref{thmjoint} (i).

\medskip

\noi (iii) Let $N/n^{2\beta} \to \infty$ and $0 < \beta < 1/2$. Then
\begin{eqnarray}\label{case3}
&N^{2-1/\beta}\widehat \gamma_{N,n} (t,s) \to_{\rm d} -(\bar V_{2\beta})^2,
\end{eqnarray}
where $\bar V_{2\beta}$ is symmetric $(2\beta)$-stable r.v. defined in Proposition \ref{propX} (ii).

\medskip

\noi (iv) Let $N/n^{2\beta} \to 0, N/n^{\beta} \to \infty$ and $3/4 < \beta < 3/2 $. Then
\begin{eqnarray} \label{case4}
&N^{1 - 3/(4\beta)} n^{1/2} \widehat \gamma_{N,n} (t,s) \ \to_{\rm d} \ W_{4\beta/3},
\end{eqnarray}
where $W_{4\beta/3}$ is a symmetric $(4\beta/3)$-stable r.v. with
characteristic function $\E \e^{\i \theta W_{4\beta/3}} = \e^{-(\sigma_0/ 2^{2\beta/3}) |\theta|^{4\beta/3}} $
and $\sigma_0 $ is the same constant as in Theorem \ref{thmjoint} (ii).

\medskip

\noi (v) Let $N/n^{2\beta} \to 0, 1/2 < \beta < 3/4 $ and $N/n^{2\beta/(4\beta -1)} \to \infty$.
Then the convergence in \eqref{case4} holds.

\medskip

\noi (vi) Let $N/n^{\beta} \to \infty, 1/2 < \beta < 3/4 $ and $N/n^{2\beta/(4\beta -1)} \to 0$.
Then the convergence in \eqref{case3} holds.

\medskip

\noi (vii) Let $N/n^{2\beta} \to 0, N/n^{\beta} \to \infty$ and $0 < \beta < 1/2 $.
Then the convergence in \eqref{case3} holds.

\medskip

\noi (viii) Let $N/n^{\beta} \to 0$ and $3/4 < \beta < 3/2 $.
Then the convergence in \eqref{case4} holds.

\medskip

\noi (ix) Let $N/n^{\beta} \to 0, 0 < \beta < 3/4 $ and $N/n^{2\beta/(5-4\beta)} \to \infty $.
Then
\begin{eqnarray}\label{case9}
&N^{2-2/\beta}\widehat \gamma_{N,n} (t,s) \to_{\rm d} -(\bar W_{\beta})^2,
\end{eqnarray}
where $\bar W_{\beta}$ is a symmetric $\beta$-stable r.v. defined in Proposition \ref{propX} (iii).

\medskip

\noi (x) Let $0 < \beta < 3/4 $ and $N/n^{2\beta/(5-4\beta)} \to 0 $.
Then the convergence in \eqref{case4} holds.

\medskip

\noi (xi) For $3/2 < \beta < 2$, let $N/n^\beta \to [0,\infty]$ and for $\beta > 2$, let $N,n \to \infty$ in arbitrary way. Then
\begin{eqnarray}\label{case11}
&N^{1/2} n^{1/2} \widehat \gamma_{N,n} (t,s) \to_{\rm d} N(0, \sigma^2),
\end{eqnarray}
where $\sigma^2 $ is given as in Theorem \ref{thmclt}.

\end{corollary}

The proof of Theorem \ref{thmjoint} in
cases (i)--(iii) is given subsections \ref{sec:muconst}-\ref{sec:muinf}.
To avoid excessive notation,
the discussion is limited to the case $(t,s) = (0,1)$
or the partial sums process $S_{N,n}(\tau) := \sum_{i=1}^N \sum_{t=1}^{\lfloor n\tau\rfloor} X_i (t) X_{i+1}(t)$. Later on
we shall extend them to general case $(t,s), s \ne 0$.

Let us give an outline of the proof of Theorem \ref{thmjoint}. Similarly to \cite{pils2014} we use the method
of characteristic function combined with `vertical' Bernstein's blocks, due to the fact that
$S_{N,n} $ is not a sum of row-independent summands as in \cite{pils2014}. Write
\begin{equation}\label{Sdecomp}
S_{N,n}(\tau) = S_{N,n;q} (\tau) + S^\dag_{N,n;q} (\tau) + S^\ddag_{N,n;q} (\tau),
\end{equation}
where the main term
\begin{eqnarray}\label{Sp}
S_{N,n;q} (\tau)&:=&\sum_{k=1}^{\tilde N_q} Y_{k,n;q} (\tau), \quad
Y_{k,n;q} (\tau) \ := \ \sum_{(k-1)q < i < kq} \sum_{t=1}^{\lfloor n\tau\rfloor} X_i (t) X_{i+1}(t), \ 1\le k \le \tilde N_q := \big\lfloor\frac{N}{q}\big\rfloor,
\end{eqnarray}
is a sum of $\tilde N_q$ `large' blocks of size $q-1 $ with
\begin{equation}\label{qNn}
q \equiv q_{N,n} \to \infty \qquad \text{as} \quad N, n \to \infty.
\end{equation}
The convergence rate of $q \in \N$ in \eqref{qNn} will be slow enough (e.g., $q = O(\log N)$)
and specified later on. The two other terms
in the decomposition \eqref{Sdecomp},
\begin{eqnarray}\label{Sdag}
&S^\dag_{N,n;q} (\tau)\ := \ \sum_{k=1}^{\tilde N_q} \sum_{t=1}^{\lfloor n\tau\rfloor} X_{kq} (t) X_{kq+1}(t), \quad
S^\ddag_{N,n;q} (\tau)\ := \ \sum_{q \tilde N_q < i \le N} \sum_{t=1}^{\lfloor n\tau\rfloor} X_{i} (t) X_{i+1}(t),
\end{eqnarray}
contain respectively $\tilde N_q = o(N)$ and $N- q\tilde N_q < q = o(N)$ row sums and will be shown to be negligible.
More precisely, we show that in each case (i)--(iii) of Theorem \ref{thmjoint},
\begin{eqnarray}
&&A^{-1}_{N,n} S_{N,n;q} (\tau)\ \to_{\rm fdd} \ {\cal S}_\beta (\tau), \label{S1} \\
&&A^{-1}_{N,n} S^\dag_{N,n;q} (\tau) = o_{\rm p}(1), \qquad A^{-1}_{N,n} S^\ddag_{N,n;q} (\tau) = o_{\rm p}(1), \label{S2}
\end{eqnarray}
where $A_{N,n} $ and ${\cal S}_\beta $ denote the normalization and the limit process, respectively, particularly,
\begin{eqnarray}
A_{N,n} \ := \ n^2 \begin{cases}
\lambda_{N,n}^{\beta}, &\lambda_\infty = \infty, \ 1 < \beta < 3/2, \\
\lambda_{N,n} (\log \lambda_{N,n})^{1/(2\beta)}, &\lambda_\infty = \infty, \ 0< \beta <1, \\
\lambda_{N,n}^{3/2}, &\lambda_\infty \in [0,\infty), \ 0< \beta < 3/2.
\end{cases}
\end{eqnarray}
Note that {\it the summands $Y_{k,n;q}, 1 \le k \le \tilde N_q $ in \eqref{Sp} are independent and identically distributed}, and the limit
$ {\cal S}_\beta (\tau)$ is infinitely divisible in cases (i)--(iii) of Theorem  \ref{thmjoint}.
Hence use of characteristic functions to prove \eqref{S1}
is natural. The proofs are limited to
one-dimensional convergence at a given $\tau >0$ since the convergence of general finite-dimensional distributions follows in a similar way.
Accordingly, the proof of \eqref{S1} for fixed $\tau >0$ reduces to
\begin{eqnarray}\label{Phiqconv}
&\Phi_{N,n;q}({\theta}) \ \rightarrow \ \Phi ({\theta}), \quad \text{as} \ N, \, n \to \infty, \, \lambda_{N,n} \to
\lambda_\infty, \quad \forall \theta \in \R,
\end{eqnarray}
where
\begin{eqnarray} \label{Phip}
&\Phi_{N,n;q}(\theta)
:= \tilde N_q \E \big[\e^{\i \theta A^{-1}_{N,n}Y_{1,n;q}(\tau)} -1 \big], \qquad
\Phi(\theta) := \log \E \e^{\i \theta {\cal S}_\beta (\tau)}.
\end{eqnarray}
To prove \eqref{Phiqconv} write
\begin{eqnarray}\label{Phip1}
&A^{-1}_{N,n}Y_{1,n;q}(\tau)
= \sum_{i=1}^{q-1} y_i(\tau), \quad \text{where} \quad y_i(\tau) := A^{-1}_{N,n}\sum_{t=1}^{\lfloor n\tau\rfloor} X_i (t) X_{i+1}(t).
\end{eqnarray}
We use the identity:
\begin{eqnarray}\label{ident}
&\prod_{1 \le i <q} (1 + w_i) - 1 = \sum_{1\le i < q} w_i + \sum_{|D| \ge 2} \prod_{i \in D} w_i,
\end{eqnarray}
where the sum $\sum_{|D| \ge 2}$ is taken over all subsets $D \subset \{1, \dots, q-1\} $ of cardinality
$|D| \ge 2 $. Applying \eqref{ident} with $w_i = \e^{\i \theta y_i(\tau)} -1 $ we obtain
\begin{eqnarray}\label{Phiqrep}
&\Phi_{N,n;q}(\theta)
:= \tilde N_q (q-1) \big[\E \e^{\i \theta y_1(\tau)} -1 \big]
+ \tilde N_q \sum_{|D| \ge 2} \E \prod_{i \in D} \big[ \e^{\i \theta y_i(\tau)} -1 \big].
\end{eqnarray}
Thus, since $\tilde N_q (q-1)/N \to 1$, \eqref{Phiqconv} follows from
\begin{eqnarray}\label{E1}
N\big[\E \e^{\i \theta y_1(\tau)} -1 \big]
&\to&\Phi (\theta), \\
N \sum_{
|D|\ge 2} \E \prod_{i \in D} \big[ \e^{\i \theta y_i(\tau)} -1 \big]
&\to&0. 
\label{E2}
\end{eqnarray}

Let us explain the main idea of the proof of \eqref{E1}. Assuming
$ \phi (x) = (1-x)^{\beta-1} $ in \eqref{a:beta}
the l.h.s.\ of \eqref{E1} can be written as
\begin{eqnarray}
N\big[\E \e^{\i \theta y_1(\tau)} -1 \big]
&=&N \int_{(0,1]^2} \E \big[\e^{\i \theta y_1(\tau)} -1
\big|a_i = 1-z_i, i=1,2 \big] (z_1 z_2)^{\beta-1} \d z_1 \d z_2 \nn \\
&=&\frac{N}{B^{2\beta}_{N,n}}\int_{(0,B_{N,n}]^2}
\E \big[\e^{\i \theta z_{N,n}(\tau;x_1,x_2)} -1 \big] (x_1 x_2)^{\beta-1} \d x_1 \d x_2,
\label{Phirep}
\end{eqnarray}
where
\begin{eqnarray} \label{zndef}
z_{N,n}(\tau;x_1,x_2)
&:=&A^{-1}_{N,n}\sum_{s_1,s_2 \in \Z} \vep_1(s_1) \vep_2(s_2)
\sum_{t=1}^{\lfloor n\tau\rfloor} \prod_{i=1}^2 \big(1 - \frac{x_i}{B_{N,n}}\big)^{t-s_i}\1(t \ge s_i) 
\end{eqnarray}
and $B_{N,n} \to \infty$ is a scaling factor of the autoregressive coefficient.
In cases (ii) and (iii) of Theorem \ref{thmjoint} (proof of
\eqref{muzero} and \eqref{interlim}) we choose this scaling factor $B_{N,n} = N^{1/(2\beta)} $ so
that $\frac{N}{B^{2\beta}_{N,n}} =1 $ and prove that the integral in
\eqref{Phirep} converges to $\int_{\R^2_+}
\E \big[\e^{\i \theta z(\tau;x_1,x_2)} -1 \big] (x_1 x_2)^{\beta-1} \d x_1 \d x_2 = \Phi(\theta)$, where
$z(\tau; x_1, x_2)$ is a random process and
$\Phi(\theta)$ is the required limit in \eqref{Phiqconv}. A similar scaling
$B_{N,n} = (N \log \lambda_{N,n})^{1/(2\beta)}$ applies in the case $\lambda_\infty = \infty$, $0< \beta < 1$
(proof of \eqref{mulog}) although in this case the factor $N/B^{2\beta}_{N,n} = 1/\log \lambda_{N,n}$ in front of the integral in \eqref{Phirep} does not trivialize and the proof of the limit in \eqref{Phiqconv} is more delicate. On the other hand, in the case of the Gaussian limit \eqref{muinf}, the choice $B_{N,n}=n$ leads to $N/B^{2\beta}_{N,n} = \lambda^{2\beta}_{N,n} \to \infty $
and \eqref{Phirep} tends to $(1/2) |\theta|^2 \int_{\R^2_+} \E z^2 (\tau; x_1, x_2) (x_1 x_2)^{\beta-1} \d x_1 \d x_2 = \Phi(\theta) $
with $z(\tau; x_1, x_2) $ defined in \eqref{zdef}
as shown in subsection \ref{sec:muinf} below.

To summarize the above discussion: in each case (i)--(iii) of Theorem \ref{thmjoint},
to prove the limit \eqref{S1} of the main term,
it suffices to verify relations \eqref{E1} and \eqref{E2}. The proof of
the first relation in \eqref{S2} is very similar to \eqref{S1} since
$S^\dag_{N,n;q} (\tau) $ is also a sum of i.i.d.\ r.v.s and the argument
of \eqref{S1} applies with small changes. The proof of the second relation
in \eqref{S2} seems even simpler. In the proofs we repeatedly use the following inequalities:
\begin{equation}\label{ineq:rho}
| \e^{\i z} - 1 | \le 2 \wedge |z|, \qquad | \e^{\i z} - 1 - \i z | \le (2 |z|) \wedge (z^2/2), \quad z \in \R.
\end{equation}


\subsection{Proof of Theorem \ref{thmjoint} (iii): case $0<\lambda_\infty < \infty$}
\label{sec:muconst}

\noi {\it Proof of \eqref{E1}.} For notational brevity, we assume $\lambda_{N,n} = \lambda_\infty = 1 $ since the general case as in \eqref{mucond} requires unsubstantial changes. Recall from \eqref{Zchf} that $\Phi(\theta) = \int_{\R^2_+} \E [\e^{\i \theta z(\tau; x_1, x_2)} -1 ] (x_1 x_2)^{\beta -1} \d x_1 \d x_2 $,
where $z(\tau; x_1, x_2)$ is the double It\^o-Wiener integral in \eqref{zdef}.
Also recall the representation \eqref{Phirep}, \eqref{zndef}, where $A_{N,n} = n^2$, $B_{N,n} = n$ and $z_{N,n}(\tau; x_1,x_2)
= Q_{12}(h_n(\cdot;\tau; x_1,x_2))$ is a quadratic form as in \eqref{Qij} with coefficients
\begin{eqnarray} \label{hn}
h_n(s_1,s_2;\tau; x_1,x_2)
&:=&n^{-2} \sum_{t=1}^{\lfloor n\tau\rfloor} \prod_{i=1}^2 \big(1 - \frac{x_i}{n}\big)^{t-s_i} \1(t \ge s_i), \qquad s_1, s_2 \in \Z.
\end{eqnarray}
By Proposition \ref{Qijconv}, with $\alpha_1 = \alpha_2 = n$,
the point-wise convergence
\begin{eqnarray}
&&\E [\e^{\i \theta z_{N,n}(\tau;x_1,x_2)} -1] = \E [\e^{\i \theta Q_{12}(h_n(\cdot;\tau; x_1,x_2))} -1]
\ \to\ \E [\e^{\i \theta z(\tau;x_1,x_2)} -1]
\label{Qnconv}
\end{eqnarray}
for any fixed $x_1, x_2 \in \R_+$ 
follows from $L_2$-convergence of the kernels:
\begin{equation} \label{hnconv1}
\|\widetilde h_n(\cdot; \tau; x_1,x_2) - h(\cdot; \tau; x_1,x_2)\| \ \to \ 0,
\end{equation}
where 
\begin{eqnarray}
\widetilde h_n(s_1,s_2; \tau; x_1,x_2)
&:=&n h_n(\lfloor ns_1 \rfloor, \lfloor ns_2 \rfloor; \tau; x_1,x_2)
\ = \ \frac{1}{n}  \sum_{t=1}^{\lfloor n\tau \rfloor} \prod_{i=1}^2 \big(1 - \frac{x_i}{n}\big)^{t-\lfloor ns_i \rfloor} \1(t \ge \lfloor ns_i \rfloor) \nn \\
&\to&\int_0^\tau \prod_{i=1}^2 \e^{-x_i(t-s_i)} \1( t> s_i) \d t =: h(s_1,s_2; \tau; x_1,x_2)\label{hnconv2}
\end{eqnarray}
point-wise for any $x_i > 0, s_i \in \R, s_i \ne 0, i=1,2, \tau > 0$ fixed. We also use the dominating
bound
\begin{equation} \label{hdomin}
|\widetilde h_n(s_1,s_2; \tau; x_1,x_2)| \le C h(s_1,s_2; 2\tau; x_1,x_2), \qquad s_1,s_2 \in \R, \qquad 0< x_1,x_2 < n,
\end{equation}
with $C >0$ independent of $s_i,x_i, i=1,2 $
which follows from the definition of $\widetilde h_n(\cdot; \tau; x_1,x_2) $ and the inequality $1-x \le \e^{-x}, x >0$. Since
$h(\cdot; 2\tau; x_1,x_2) \in L^2(\R^2) $, \eqref{hnconv2}, \eqref{hdomin} and the dominated convergence theorem imply
\eqref{hnconv1} and \eqref{Qnconv}.

It remains to show the convergence of the corresponding integrals, viz.,
\begin{equation}\label{intconv}
\int_{(0,n]^2} \E [\e^{\i \theta z_{N,n}(\tau;x_1,x_2)} -1 ] (x_1 x_2)^{\beta-1} \d x_1 \d x_2
\to \int_{\R^2_+} \E [\e^{\i \theta z(\tau;x_1,x_2)} -1 ] (x_1 x_2)^{\beta-1} \d x_1 \d x_2 = \Phi(\theta).
\end{equation}
From \eqref{Phirep} and $\E z_{N,n}(\tau;x_1,x_2) = 0 $ we obtain
\begin{eqnarray}\label{Ezbdd}
\big|\E \big[\e^{\i \theta z_{N,n}(\tau;x_1,x_2)} -1 \big] \big|
&\le&C\begin{cases}
1, &x_1 x_2 (x_1 + x_2) \le 1, \\
\E z^2_{N,n}(\tau;x_1,x_2), &x_1 x_2(x_1+x_2) > 1,
\end{cases}
\end{eqnarray}
where
\begin{eqnarray}
\E z^2_{N,n}(\tau;x_1,x_2)&=&A^{-2}_{N,n} \E \Big[ \Big(\sum_{t=1}^{\lfloor n\tau\rfloor} Y_{12}(t)\Big)^2 |a_i = 1 - \frac{x_i}{B_{N,n}}, i=1,2 \Big]\nn \\
&=&n^{-4} \E \Big[ \Big(\sum_{t=1}^{\lfloor n\tau\rfloor} Y_{12}(t)\Big)^2|a_i = 1 - \frac{x_i}{n}, i=1,2\Big]\nn \\
&\le&\frac{C}{n^3 (x_1/n) (x_2/n)}
 \min \big\{n, \frac{1}{ (x_1 + x_2)/n}\big\} \ = \
\frac{C}{x_1 x_2} \min \big\{1, \frac{1}{x_1 + x_2}\big\},
 \label{Ez2}
\end{eqnarray}
see \eqref{zndef} and the bound in \eqref{var1}. In view of inequality \eqref{Aineq2},
the dominated convergence theorem applies, proving \eqref{intconv} and \eqref{E1}.
%
%
%

\medskip

\noi {\it Proof of \eqref{E2}}. Choose $q = q_{N,n} = \lfloor\log n\rfloor $. Let $J_q(\theta) $ denote the l.h.s.\ of
\eqref{E2}. Using the identity $\sum_{D \subset \{1, \dots, q-1\}: |D| \ge 2} \prod_{i \in D} w_i
= \sum_{1 \le i < j < q} w_i w_j \prod_{i < k < j} (1+ w_k) $ with $w_i = \e^{\i \theta y_i(\tau)} -1 $, see \eqref{ident},
we can rewrite $J_q(\theta) = \sum_{1 \le i < j < q} T_{ij}(\theta)$, where
\begin{eqnarray}\label{Tijdef}
T_{ij}(\theta) &:=& N \E \Big[  (\e^{\i \theta y_i(\tau)} -1)(\e^{\i \theta y_j(\tau)} -1) \exp \Big\{\i \theta \sum\nolimits_{i<k < j} y_k(\tau)\Big\} ( \1(a_i < a_{j+1}) + \1(a_i > a_{j+1}) ) \Big]\\
&=& T'_{ij} (\theta) + T''_{ij}(\theta).\nn
\end{eqnarray}
Since $|J_q(\theta)| \le q^2 \max_{1 \le i < j < q} |T_{ij}(\theta)| \le (\log n)^2  \max_{1 \le i < j < q} |T_{ij}(\theta)|$,
\eqref{E2} follows from
\begin{equation}\label{Tij}
 |T_{ij}(\theta)| \le C n^{-\delta}, \qquad \forall \ 1 \le i < j,
\end{equation}
with $C, \delta >0$ independent of $n $.
Using $\E [y_i(\tau)|a_k, \vep_j(k),
k, j \in \Z, j \ne i ] = 0$ and \eqref{Ez2} we obtain 
\begin{eqnarray}\label{ineq:Tijbdd}
|T'_{ij}(\theta)| &\le&CN \E \big[ \min \big\{1, \E [y^2_i(\tau)|a_k, k \in \Z]\big\}
 \1(a_i < a_{j+1}) \big] \\
&\le&C n^{-\beta} \int_{(0,n]^3} \min \Big\{1, \frac{1}{x_i x_{i+1} (x_i + x_{i+1})} \Big\} (x_i x_{i+1} x_{j+1})^{\beta-1} \1(x_{j+1} < x_i) \d x_i \d x_{i+1} \d x_{j+1} \nn \\
&=&Cn^{-\beta} \int_{(0,n]^2} \min \Big\{1, \frac{1}{x_1 x_2 (x_1 + x_2)} \Big\} x^{2\beta -1}_1 x_2^{\beta-1} \d x_1 \d x_2 \le C n^{-\beta} (T'_n + T''_n), \nn
\end{eqnarray}
where $T'_n := \int_{0 < x_1 < x_2 < n} \min \big\{1, \frac{1}{x_1 x^2_2} \big\} x^{2\beta -1}_1 x_2^{\beta-1} \d x_1 \d x_2$,
$T''_n := \int_{0 < x_2 < x_1 < n} \min \big\{1, \frac{1}{x^2_1 x_2} \big\}x^{2\beta -1}_1 x_2^{\beta-1} \d x_1 \d x_2 $.
Next,
\begin{eqnarray*}
T'_n&\le&\int_0^1 x_1^{2\beta -1} \d x_1 \Big[\int_{x_1}^{1/\sqrt{x_1}} x_2^{\beta-1} \d x_2
+ x_1^{-1}\int_{1/\sqrt{x_1}}^n x_2^{\beta-3} \d x_2 \Big] +
\int_1^n x_1^{2\beta -2} \d x_1 \int_{x_1}^n x_2^{\beta -3} \d x_2 \\
&\le&C \Big[\int_0^1 x_1^{3\beta/2 -1} \d x_1 + \int_1^n x_1^{3\beta - 4} \d x_1 \Big]
\ \le \ C n^{3(\beta -1)\vee 0} (1 + \1(\beta =1)\log n ).
\end{eqnarray*}
Similarly,
\begin{eqnarray*}
T''_n &=&\int_0^1 x_1^{2\beta-1} \d x_1 \int_0^{x_1} x_2^{\beta-1} \d x_2 + \int_1^n x_1^{2\beta-1} \d x_1 \int_0^{x_1^{-2}} x_2^{\beta-1} \d x_2 + \int_1^n x_1^{2\beta-3} \d x_1 \int_{x_1^{-2}}^{x_1} x_2^{\beta-2} \d x_2\\
&\le& C \big( (\log n) \1 (\beta < 1) + (\log n)^2 \1 (\beta=1) + n^{3(\beta-1)} \1 (\beta > 1) \big).
\end{eqnarray*}
Whence, the bound in  \eqref{Tij} follows for $T'_{ij}(\theta)$ with any
$0< \delta < \beta \wedge (3-2\beta)$, for $0< \beta < 3/2 $.
Since $ |T''_{ij}(\theta)|
\le CN \E [ \min \{1, \E [y^2_j(\tau)|a_k, k \in \Z] \} \1(a_{j+1} < a_{i}) ]$
can be symmetrically handled, this proves
\eqref{Tij} and \eqref{E2}.

\medskip

\noi \noi {\it Proof of \eqref{S2}.} Since $A^{-1}_{N,n} S^\dag_{N,n;q} (\tau) = \sum_{k=1}^{\tilde N_q} y_{kq}(\tau) $ is a
sum of $\tilde N_q $ i.i.d.\ r.v.s $y_{kq}(\tau), k=1, \dots, \tilde N_q $, so the first relation
in \eqref{S2} follows from
\begin{equation}\label{ES1}
\tilde N_q \E [\e^{\i \theta y_1(\tau)} -1] \ \to \ 0, \qquad \forall \ \theta \in \R.
\end{equation}
Clearly, \eqref{ES1} is a direct consequence of \eqref{E1} and the fact that $\tilde N_q/N \to 0$.

Consider the second relation in \eqref{S2}. Let $L_q := N - q \tilde N_q $ be the number of summands
in $S^\ddag_{N,n;q} (\tau)$. Then $ A^{-1}_{N,n} S^\ddag_{N,n;q} (\tau) =_{\rm fdd}
\sum_{i=1}^{L_q} y_i(\tau) $ and
\begin{eqnarray}\label{ES2}
\E \e^{\i \theta A^{-1}_{N,n} S^\ddag_{N,n;q} (\tau)} -1
&=&L_q \E [\e^{\i \theta y_1(\tau)} -1] + \sum_{|D|\ge 2} \E \prod_{i \in D} \big[ \e^{\i \theta y_i(\tau)} -1 \big],
\end{eqnarray}
where the last sum is taken over all $D \subset \{1, \dots, L_q\}, |D| \ge 2 $. Since $L_q < q = o(N)$
from \eqref{E1}, \eqref{E2} we infer that the r.h.s.\ of \eqref{ES2} vanishes, proving
\eqref{S2}, and thus completing the proof of  Theorem \ref{thmjoint}, case (iii).


\subsection{Proof of Theorem \ref{thmjoint} (ii): case $\lambda_\infty = 0$, or $N = o(n^{2\beta})$.}
\label{sec:mu0}

\noi {\it Proof of \eqref{E1}.}
Note the log-ch.f. of the r.h.s.\ in \eqref{muzero} can be written as
\begin{eqnarray}
\Phi(\theta)&=&\log \E \e^{\i \theta {\cal A}^{1/2} B(\tau)}\ = \
\log \E \e^{ - (\theta^2 \tau/2) {\cal A} }
= - \sigma_0 (\theta^2 \tau / 2)^{2\beta/3} \nn \\
&=&- \mbox{$\psi(1)^2 \int_{\R_+^2} \Big( 1 - \exp \Big\{ - \frac{\theta^2 \tau}{4 x_1 x_2 (x_1 + x_2)}\Big\} \Big) (x_1 x_2)^{\beta-1}
\d x_1 \d x_2$} \label{stableint}
\end{eqnarray}
with $\sigma_0 >0 $ given by the integral
\begin{eqnarray}\label{stableint1}
&\sigma_0 \
:= \ \psi (1)^2 2^{-2\beta/3} \int_{\R_+^2} \Big( 1 - \exp \Big\{ - \frac{1}{x_1 x_2 (x_1 + x_2)}\Big\} \Big) (x_1 x_2)^{\beta-1} \d x_1 \d x_2.
\end{eqnarray}
Relation \eqref{stableint} follows by change of variable $x_i \to (\theta^2 \tau/4) ^{1/3} x_i, i=1,2 $. The convergence
of the integral in \eqref{stableint1} follows from \eqref{Aineq2}. The explicit
value of $\sigma_0$ in \eqref{stableint1} is given in Proposition \ref{propinter} (v) and
 computed in the Appendix.
Recall the representation in \eqref{Phirep} where $B_{N,n} = N^{1/(2\beta)}, N/B^{2\beta}_{N,n} = 1 $ and
\begin{eqnarray}\label{zNn2}
z_{N,n}(\tau;x_1,x_2)
&:=&N^{-3/(4\beta)} n^{-1/2} \sum_{s_1,s_2 \in \Z} \vep_1(s_1) \vep_2(s_2)
\sum_{t=1}^{\lfloor n\tau\rfloor} \prod_{i=1}^2 \Big(1 - \frac{x_i}{N^{1/(2\beta)}}\Big)^{t-s_i}\1(t \ge s_i). 
\end{eqnarray}
Let us prove the (conditional) CLT:
\begin{equation}\label{zlim}
z_{N,n}(\tau;x_1,x_2) \ \to_{\rm fdd} \ (2 x_1 x_2 (x_1 + x_2))^{-1/2} B(\tau),
\end{equation}
implying the point-wise convergence
\begin{eqnarray}
&\E[1- \e^{\i \theta z_{N,n}(\tau;x_1,x_2)}]
\ \to\ 1 - \e^{- \theta^2 \tau/(4 x_1 x_2 (x_1 +x_2))}
\label{subQnconv}
\end{eqnarray}
of the integrands in \eqref{Phirep} and \eqref{stableint1},
for any fixed $(x_1,x_2) \in \R^2_+$. As in the rest of the paper, we restrict the proof of \eqref{zlim}
to one-dimensional convergence, and set $\tau =1 $ for concreteness. Split \eqref{zNn2} as
$z_{N,n}(1;x_1,x_2) = z^+_{N,n}(x_1,x_2) + z^-_{N,n}(x_1,x_2)$, where
$z^+_{N,n}(x_1,x_2)
:= N^{-3/4\beta} n^{-1/2} \sum_{s_1,s_2 =1}^n \vep_1(s_1) \vep_2(s_2) \cdots $ corresponds to
the sum over $1 \le s_1, s_2 \le n $ alone. Thus, we shall prove that
\begin{equation} \label{zzlim}
z^-_{N,n}(x_1,x_2) = o_{\rm p}(1) \quad \text{and} \quad z^+_{N,n}(x_1,x_2) \ \to_{\rm d} \ N\Big(0, \frac{1}{2 x_1 x_2 (x_1 + x_2)}\Big).
\end{equation}
Arguing as in the proof of \eqref{var1} it is easy to show that
\begin{eqnarray*}
\E (z^-_{N,n}(x_1,x_2))^2&\le&\frac{C}{N^{3/(2\beta)}n} \Big(\frac{x_1+ x_2}{N^{1/(2\beta)}} \Big)^{-2}
\Big\{ \Big(\frac{x_1}{N^{1/(2\beta)}}\Big)^{-2} + \Big(\frac{x_2}{N^{1/(2\beta)}}\Big)^{-2}
+ \Big(\frac{x_1}{N^{1/(2\beta)}}\Big)^{-1} \Big(\frac{x_2}{N^{1/(2\beta)}}\Big)^{-1} \Big\} \nn \\
&=&C \lambda_{N,n} (x_1 + x_2)^{-2} \big\{x_1^{-2} + x_2^{-2} + (x_1x_2)^{-1}\big\}, \label{z2bdd}
\end{eqnarray*}
where $\lambda_{N,n} \to 0 $, implying the first relation in \eqref{zzlim}.
To prove the second relation in \eqref{zzlim}
we use the martingale CLT in Hall and Heyde \cite{hall1980}. (The same approach is used to prove CLT
for quadratic forms in \cite{bhan2007}.)
Towards this aim,  
write $z^+_{N,n}(x_1,x_2)$ as a sum of zero-mean square-integrable martingale difference array
\begin{eqnarray*}
&z^+_{N,n}(x_1,x_2)=\sum_{k=1}^n Z_k, \quad
Z_k := \vep_1(k)\sum_{s=1}^{k-1} f(k,s)\, \vep_2 (s) + \vep_2(k) \sum_{s=1}^{k-1} f(s,k)\, \vep_1(s) + \vep_1 (k)\vep_2(k) f(k,k)
\end{eqnarray*}
with respect to the filtration ${\cal F}_k$ generated by $\{\vep_i (s), \, 1\le s \le k, i=1,2\}$, $0\le k \le n$, where
\begin{eqnarray*}
f(s_1,s_2) &:=& N^{-3/(4\beta)} n^{-1/2}\sum_{t=1}^n \prod_{i=1}^2 \big(1 - \frac{x_i}{N^{1/(2\beta)}}\big)^{t-s_i}\1(t \ge s_i), \quad
1 \le s_1,s_2 \le n.
\end{eqnarray*}
Accordingly, the second convergence in \eqref{zzlim} follows from
\begin{equation}\label{a:CLT}
\sum_{k=1}^n \E [Z^2_{k} | {\cal F}_{k-1}] \to_{\rm p} \frac{1}{2 x_1 x_2 (x_1+x_2)} \quad \text{and} \quad
\sum_{k=1}^n \E [Z^2_{k} \1 (|Z_{k}|> \epsilon )] \to 0 \quad \text{for any } \epsilon > 0.
\end{equation}
Note the conditional variance
$v^2_k := \E [Z_{k}^2 | {\cal F}_{k-1}]=\big(\sum_{s=1}^{k-1} f(k,s) \vep_2(s)\big)^2 + \big(\sum_{s=1}^{k-1} f(s,k)
\vep_1(s)\big)^2 + f^2(k,k)$, where
\begin{equation} \label{EZlim}
\sum_{k=1}^n \E Z^2_k \ = \ \sum_{k=1}^n \E v^2_k \ = \
\sum_{s_1,s_2=1}^n f^2(s_1,s_2)\ = \ \E (z^+_{N,n}(x_1,x_2))^2 \ \to \ \frac{1}{2 x_1 x_2 (x_1+x_2)}
\end{equation}
is a direct consequence of the asymptotics in \eqref{var2},
where $a_i=1-x_1/N^{1/(2 \beta)}, a_j = 1 - x_2/N^{1/(2\beta)}$.
 Therefore the first relation in \eqref{a:CLT} follows from
\eqref{EZlim} and
\begin{equation} \label{limR}
R_n := \sum_{k=1}^n (v^2_k - \E v^2_k) = o_{\rm p}(1).
\end{equation}
To show \eqref{limR} we split
$R_n = R'_n + R''_n $ into the sum of `diagonal' and `off-diagonal' parts, viz.,
\begin{eqnarray*}
R'_n&:=&\sum_{i=1}^2 \sum_{1\le s < n} c_i(s) (\vep_i^2(s)-1), \qquad R''_n\ :=\ \sum_{i=1}^2 \sum_{1\le s_1, s_2 < n, s_1 \ne s_2} c_i(s_1,s_2) \vep_i(s_1) \vep_i (s_2),
\end{eqnarray*}
where $ c_{1}(s) := \sum_{s < k \le n} f^2(s,k)$, $c_{2}(s):=\sum_{s < k \le n} f^2(k,s)$, $
c_1(s_1,s_2) := \sum_{s_1 \vee s_2 < k \le n} f(s_1,k) f(s_2,k)$,
$c_2(s_1,s_2) := \sum_{s_1 \vee s_2 < k \le n} f(k,s_1) f(k,s_2)$.
Using the elementary bound for $1 \le s_1, s_2 \le n$: $\sum_{t=1}^n \prod_{i=1}^2 a_i^{t-s_i} \1(t \ge s_i)
\le \big(a_2^{s_1 - s_2} \1(1 \le s_2 \le s_1) +
a_1^{s_2 - s_1} \1(1 \le s_1 \le s_2) \big) S(a_1,a_2)$,
$S(a_1,a_2) := \sum_{t=0}^\infty (a_1 a_2)^t = (1-a_1 a_2)^{-1} \le 2(2 - a_1 - a_2)^{-1}$,
we obtain
\begin{eqnarray} \label{cbdd}
&&|c_i(s)| \le Cn^{-1} x_i^{-1} (x_1 + x_2)^{-2}, \qquad
\sum_{s_1,s_2=1}^n c^2_i(s_1,s_2) \le C \lambda_{N,n} x_i^{-3} (x_1 + x_2)^{-4}, \qquad i=1,2.
\end{eqnarray}
By \eqref{cbdd}, 
for $1 < p < 2 $ and $x_1, x_2 > 0$ fixed
\begin{eqnarray}
&&\E |R'_n|^p \le
C\sum_{i=1}^2 \sum_{s=1}^{n-1} |c_i(s)|^p \ \le \ C
n^{-(p-1)} = o(1), \label{R1bdd}\\
&&\E |R''_n|^2\ \le\ \sum_{i=1}^2 \sum_{s_1,s_2=1}^n c^2_i(s_1,s_2) \le C \lambda_{N,n} = o(1),
\end{eqnarray}
proving \eqref{limR} and the first relation in \eqref{a:CLT}. The proof of the second relation in \eqref{a:CLT}
is similar since it reduces to $T_n := \sum_{k=1}^n \E [ |Z_{k}|^{2p}] = o(1) $ for the same $1< p \le 2 $, where
$\E |Z_k|^{2p} \le C \big(\E |\sum_{s=1}^{k-1} f(k,s)\, \vep_2 (s)|^{2p} + \E |\sum_{s=1}^{k-1} f(s,k)\, \vep_1(s)|^{2p} + |f(k,k)|^{2p}\big)
\le C \big( (\sum_{s=1}^{k-1} f^2(k,s))^p + (\sum_{s=1}^{k-1} f^2(s,k))^p + |f(k,k)|^{2p} \big)$ by Rosenthal's inequality, see e.g.
(\cite{gir2012}, Lemma 2.5.2), and the sum $T_n = O(n^{-(p-1)}) = o(1)$ similarly to \eqref{R1bdd}.
This proves  \eqref{a:CLT}, \eqref{zzlim}, and the pointwise convergence in
\eqref{subQnconv}.

Now we return to the proof of \eqref{E1}, whose both sides are written
as respective integrals \eqref{Phirep} and \eqref{stableint}. Due to the convergence
of the integrands (see \eqref{subQnconv}), it suffices to justify the passage to the limit using a dominated convergence theorem argument.
The dominating function independent of $N, n$ is obtained from
\eqref{Phirep} and $\E z_{N,n}(\tau;x_1,x_2) = 0 $ and
from \eqref{Ezbdd}, \eqref{Ez2}, \eqref{Aineq2} similarly as in the case $\lambda_\infty \in (0,\infty)$ above.
This proves \eqref{E1}.

\smallskip

{\it Proofs of \eqref{E2} and \eqref{S2}} are completely analogous to those in the case $\lambda_\infty \in (0,\infty)$
except that we now choose $q = \lfloor \log N \rfloor$ and replace $n$ in \eqref{Tij} and elsewhere in the proof of \eqref{E2} and \eqref{S2}, case $\lambda_\infty \in (0,\infty)$,
by $N^{1/2\beta} $.
This ends the proof of  Theorem \ref{thmjoint}, case (ii).

\subsection{Proof of Theorem \ref{thmjoint} (i): case $\lambda_\infty = \infty$, or $n = o(N^{1/(2\beta)})$} \label{sec:muinf}

\noi {\it Case $1 < \beta < 3/2$. Proof of \eqref{E1}.} In this case, $\Phi(\theta) := - \sigma^2_\infty \tau^{2(2-\beta)} \theta^2/2,
B_{N,n} = n $ and $A_{N,n} = n^2 \lambda^\beta_{N,n} = n^{2-\beta} N^{1/2}$.
Rewrite the l.h.s. of \eqref{E1} as
\begin{eqnarray}\label{GNn}
&N\big[\E \e^{\i \theta y_1(\tau)} -1 \big]
= \int_{[0,n)^2} \E \Lambda_{N,n}(\theta; \tau; x_1,x_2) (x_1 x_2)^{\beta-1} \d x_1 \d x_2, \qquad \text{where} \\
&\Lambda_{N,n}(\theta; \tau; x_1,x_2) := \lambda^{2\beta}_{N,n}
\big[\e^{\i \theta \lambda^{-\beta}_{N,n}\tilde z_{N,n}(\tau;x_1,x_2)} -1 - \i \theta \lambda^{-\beta}_{N,n} \tilde z_{N,n}(\tau;x_1,x_2) \big] \nn 
\label{Phirep3}
\end{eqnarray}
and where $\tilde z_{N,n}(\tau;x_1,x_2) $ is defined as in \eqref{zndef} with
$A_{N,n}$ replaced by $\tilde A_{N,n} := n^2 = A_{N,n}/\lambda_{N,n}^\beta. $
As shown in the proof of Case (iii) (the `intermediate limit'), for any $x_1, x_2 > 0$
\begin{equation} \label{tizconv}
\tilde z_{N,n}(\tau;x_1,x_2) \to_{\rm d} z(\tau;x_1,x_2) \quad \text{and} \quad
\E \tilde z^2_{N,n}(\tau;x_1,x_2) \to \E z^2(\tau;x_1,x_2),
\end{equation}
see \eqref{Qnconv}, where $z(\tau;x_1,x_2)$ is defined in \eqref{zdef} and the last
expectation in \eqref{tizconv} is given in \eqref{z2}.
Then using Skorohod's representation we extend \eqref{tizconv} to $\tilde z_{N,n}(\tau;x_1,x_2) \to z(\tau;x_1,x_2)$
a.s. implying also $\Lambda_{N,n}(\tau;x_1,x_2) \to - (\theta^2/2) z^2(\tau;x_1,x_2) $
a.s. Since $|\Lambda_{N,n}(\theta; \tau; x_1, x_2)| \le C \tilde z^2_{N,n} (\tau; x_1, x_2) $
and \eqref{tizconv} holds, by Pratt's lemma we obtain
\begin{eqnarray} \label{GNnconv}
&G_{N,n}(\tau; x_1,x_2)\ \to \ -\frac12\,\E z^2(\tau;x_1,x_2), \qquad \forall \, (x_1,x_2) \in \R^2_+.
\end{eqnarray}
Relation \eqref{E1} follows from
\eqref{GNn}, \eqref{GNnconv} and the dominated convergence theorem,
using the dominating bound
\begin{eqnarray} \label{GNndom}
&|G_{N,n}(\tau; x_1,x_2)| \ \le \ C \E \tilde z^2_{N,n}(\tau;x_1,x_2)\ \le \ \frac{C}{x_1 x_2} \min \big\{ 1, \frac{1}{x_1+x_2} \big\} \ =: \ \bar G(x_1,x_2),
\end{eqnarray}
see \eqref{Ez2}, and integrability of $\bar G$, see \eqref{Ez5}.

\smallskip

\noi {\it Proof of \eqref{E2}} is similar to that in case (iii) $0<\lambda_\infty < \infty $ above with $q = \lfloor\log n\rfloor$. It suffices
to check the bound \eqref{Tij} for $T_{ij}(\theta) = T'_{ij} (\theta) + T''_{ij} (\theta)$ given in \eqref{Tijdef}. By the same argument as in \eqref{ineq:Tijbdd}, we obtain $|T'_{ij}(\theta)| \le C N \E [ y^2_i (\tau) \1 (a_i < a_{j+1}) ]$. The bound on $\E \tilde z_{N,n}^2 (\tau; x_1, x_2)$ in \eqref{GNndom} further implies
$$
|T'_{ij}(\theta)| \le C n^{-\beta} \int_{(0,n]^3} \frac{1}{x_1 x_2} \min \Big\{ 1, \frac{1}{x_1+x_2} \Big\} (x_1 x_2 x_3)^{\beta-1} \1 (x_3 < x_1 ) \d x_1 \d x_2 \d x_3 \le C n^{-\beta} (T'_n + T''_n),
$$
where
$$
T'_n := \int_0^n \min\Big\{ 1, \frac{1}{x_1} \Big\} x_1^{2\beta-2} \d x_1 \int_0^{x_1} x_2^{\beta-2} \d x_2
= C \Big( \int_0^1 x_1^{3\beta-3} \d x_1 + \int_1^n x_1^{3\beta-4} \d x_1 \Big) \le C n^{3\beta-3}
$$
and
$$
T''_n := \int_0^n \min \Big\{1, \frac{1}{x_2}\Big\} x_2^{\beta-2} \d x_2 \int_0^{x_2} x_1^{2\beta-2} \d x_1\\
= C \Big( \int_0^1 x_2^{3\beta-3} \d x_2 + \int_1^n x_2^{3\beta-4} \d x_2 \Big) \le C n^{3\beta-3}.
$$
Then $|T''_{ij} (\theta)| \le C N \E [ y_j^2 (\tau) \1 (a_i > a_{j+1}) ]$ can be handled in the same way.
Whence, the bound in \eqref{Tij} follows
with any $0< \delta < 3-2\beta$, for $1< \beta < 3/2 $.
This proves \eqref{E2}. Proof of \eqref{S2} using $\tilde N_q/N \to 0$ and $L_q = N - q \tilde N_q < q = o(N)$
is completely analogous to that in case (iii) $0<\lambda_\infty < \infty $.
This completes the proof of  Theorem \ref{thmjoint}, case (i) for $1< \beta < 3/2 $.

\medskip

\noi {\it Case $0 < \beta < 1$. Proof of \eqref{E1}.} In the rest of this proof,
write $\lambda \equiv \lambda_{N,n} = N^{1/(2\beta)}/n \to \infty $ for brevity.
Also denote $\lambda' := \lambda (\log \lambda)^{1/2\beta}$,
$\log \lambda' /\log \lambda \to 1. $ Let $B_{N,n} := \lambda' n$, then
\begin{eqnarray} \label{zndef31}
z_{N,n}(\tau;x_1,x_2)
&:=&\frac{1}{\lambda' n^2}\sum_{s_1,s_2 \in \Z} \vep_1(s_1) \vep_2(s_2)
\sum_{t=1}^{\lfloor n\tau\rfloor} \prod_{i=1}^2 \big(1 - \frac{x_i}{\lambda' n}\big)^{t-s_i}\1(t \ge s_i). 
\end{eqnarray}
Split the r.h.s.\ of \eqref{E1} as follows:
\begin{eqnarray*}
N \E \big[\e^{\i \theta y_1(\tau)} -1\big]
&=&\frac{1}{\log \lambda}
\int_{(0, \lambda' n]^2} \big(\1 (1< x_1 + x_2 < \lambda) + \1(x_1+x_2 > \lambda) + \1(x_1 + x_2 < 1)\big) \\
&&\times \E \big[\e^{\i \theta z_{N,n}(\tau;x_1,x_2)} - 1\big]
(x_1 x_2)^{\beta -1} \d x_1 \d x_2 \ =: \ \sum_{i=1}^3 L_i.
\end{eqnarray*}
Here, $L_1 $ is the main term and $L_i, i=2,3 $ are remainders.
Indeed, $|L_3| = O(1/\log \lambda) = o(1)$. To estimate $L_2$ we need the bound
\begin{eqnarray}
\E z^2_{N,n}(\tau;x_1,x_2)
\ \le \  \frac{C}{x_1 x_2} \min \Big\{1, \frac{\lambda'}{x_1 + x_2}\Big\},
\label{Qnewbdd}
\end{eqnarray}
which follows from \eqref{var1} similarly to \eqref{Ez2}.
Using \eqref{Qnewbdd} we obtain
\begin{eqnarray}\label{L2}
|L_2|
&\le&
\frac{C}{\log \lambda}
\int_{x_1 + x_2 > \lambda}
\min \Big\{1, \frac{\lambda'}{x_1 x_2(x_1 + x_2)}\Big\}
(x_1 x_2)^{\beta -1} \d x_1 \d x_2 \ = \ \frac{C}{\log \lambda}(J'_\lambda + J''_\lambda),
\end{eqnarray}
where, by change of variables: $x_1 + x_2 = y, x_1 = yz $,
\begin{eqnarray*}
J'_\lambda
&:=&\int_{x_1 + x_2 > \lambda} \1(x_1 x_2(x_1+x_2) < \lambda') (x_1 x_2)^{\beta -1} \d x_1 \d x_2 \\
&=&\int_\lambda^\infty \int_0^1
\1 (y^3 z(1-z) < \lambda') y^{2\beta -1} (z(1-z))^{\beta -1} \d z \d y \\
&\le&C\int_{\lambda }^\infty y^{2\beta -1} \d y
\int_0^{1/2} z^{\beta -1} \1 (y^3 z < 2 \lambda') \d z \
\le\ C (\lambda')^\beta \int_{\lambda}^\infty y^{-\beta -1} \d y \ = \ C (\log \lambda)^{1/2}
\end{eqnarray*}
since $0< \beta < 1 $. Similarly,
\begin{eqnarray*}
J''_\lambda
&:=&\lambda'
\int_{x_1 + x_2 > \lambda} \1(x_1 x_2(x_1+x_2) > \lambda')
(x_1 + x_2)^{-1} (x_1 x_2)^{\beta -2} \d x_1 \d x_2 \\
&\le&C\lambda' \int_{\lambda}^\infty y^{2\beta -4} \d y
\int_0^{1/2} z^{\beta -2} \1 (y^3 z > \lambda' ) \d z
\ \le \ C (\log \lambda)^{1/2}.
\end{eqnarray*}
This proves $|L_2| = O(1/\log \lambda) = o(1)$.

Consider the main term $L_1$. Although $\E \e^{\i \theta z_{N,n}(\tau; x_1,x_2)} $ and hence the integrand in $L_1$
point-wise converge for any $(x_1,x_2) \in \R^2_+$,
see below, this fact is not very useful since the contribution to the limit of $L_1$ from
bounded $x_i$'s is negligible due to the presence of the factor $ 1/\log \lambda \to 0 $ in front of this integral. It turns out 
that
the main (non-negligible) contribution to this integral comes from unbounded $x_1, x_2$ with
$x_1/x_2 + x_2/x_1 \to \infty $ and
$x_1 x_2 \to z \in \R_+. $
To see this,
by change of variables $y = x_1 + x_2, x_1 = yw $ and then $w = z/y^2$
we rewrite
\begin{eqnarray} \label{L1}
L_1
&=&\frac{1}{\log \lambda} \int_1^{\lambda} V_{N,n} (\theta; y) \frac{\d y}{y},
\end{eqnarray}
where
\begin{eqnarray}\label{Vint}
V_{N,n}(\theta; y)&:=&2\int_0^{y^2/2} 
\E \Big[ \exp \Big\{\i \theta z_{N,n}(\tau; \frac{z}{y}, y\Big(1 - \frac{z}{y^2}\Big) \Big\} - 1 \Big]
z^{\beta -1} \Big(1 - \frac{z}{y^2}\Big)^{\beta -1} \d z.
\end{eqnarray}
In view of $L_i = o(1), i=2,3 $ relation \eqref{E1} follows from representation \eqref{L1} and
the existence of the limit:
\begin{equation} \label{Vlim}
\lim_{y \to \infty, y = O(\lambda)} V_{N,n} (\theta; y) = V(\theta) := - k_\infty |\theta|^{2\theta} \tau^{2\beta},
\end{equation}
where the constant $k_\infty >0$ is defined below in \eqref{cbeta}. More precisely, \eqref{Vlim} says that for any $\epsilon >0 $ there exists $K >0$ such that for any
$N, n, y \ge K $ satisfying $ y \le \lambda, \lambda \ge K $
\begin{equation} \label{VVlim}
|V_{N,n} (\theta; y) - V(\theta)| < \epsilon.
\end{equation}
To show that \eqref{VVlim} implies $L_1 \to V(\theta) $ it suffices to
split $L_1 - V(\theta) = (\log \lambda)^{-1} \int_K^\lambda (V_{N,n}(\theta;y) - V(\theta)) \frac{\d y}{y} +
(\log \lambda)^{-1} \int_1^K (V_{N,n}(\theta;y) -V(\theta)) \frac{\d y}{y}$
and use \eqref{VVlim} together with the fact that $|V_{N,n}(\theta;y)|\le C$ is bounded uniformly in $N, n, y $.

To prove \eqref{VVlim}, rewrite $V(\theta) $ of \eqref{Vlim} as the integral
\begin{equation} \label{VVint}
V(\theta) = 2\int_0^{\infty} z^{\beta -1} \E (\e^{\i \theta \tau Z_1 Z_2/(2\sqrt{z})} -1) \d z
= - 2\E \int_0^{\infty} z^{\beta -1} (1- \e^{-\theta^2 \tau^2 Z^2_1/(8z)}) \d z
= - k_\infty |\theta|^{2\beta} \tau^{2\beta}
\end{equation}
with $Z_1, Z_2 \sim N(0,1)$ independent normals and
\begin{eqnarray} \label{cbeta}
k_\infty &=& \textstyle 2\E \int_0^\infty z^{\beta -1} (1- \e^{- Z^2_1/(8z)}) \d z =
2^{1-3\beta} \E |Z_1|^{2\beta} \int_0^{\infty} z^{\beta -1} (1- \e^{-1/z}) \d z\nn \\
&=& \textstyle 2^{1-2\beta} \Gamma(\beta + \frac{1}{2})
\Gamma (1-\beta)/ (\sqrt{\pi}\beta).
\end{eqnarray}
Let $\Lambda_{N,n}(z; y) := \E \big[ \exp \big\{\i \theta z_{N,n}(\tau; \frac{z}{y}, y\big(1 - \frac{z}{y^2}\big) \big\} - 1 \big]$, $\Lambda (z) := \E [\e^{\i \theta \tau Z_1 Z_2/(2\sqrt{z})} -1]$ denote the corresponding
expectations in \eqref{Vint}, \eqref{VVint}. Clearly,
\eqref{VVlim} follows from
\begin{eqnarray}\label{Llim}
&\lim_{y \to \infty, y = O(\lambda)} \Lambda_{N,n} (z; y) = \Lambda(z), \qquad \forall \ z >0,
\end{eqnarray}
and
\begin{eqnarray}\label{Lbdd}
&|\Lambda_{N,n} (z; y)| \le C(1 \wedge (1/z)), \qquad \forall \ 0< y < \lambda, \ 0< z < y^2/2.
\end{eqnarray}
The dominating bound in \eqref{Lbdd} is a consequence of \eqref{Qnewbdd}. To show
\eqref{Llim} use Proposition \ref{Qijconv} by writing $z_{N,n}(\tau; z/y, y'), y' := y(1 - z/y^2) $ in \eqref{Vint} as
the quadratic form: $z_{N,n}(\tau; z/y, y') = Q_{12} (h_{\alpha_1,\alpha_2}(\cdot; \tau;z)) $ with
\begin{eqnarray} \label{halpha}
h_{\alpha_1,\alpha_2}(s_1,s_2;\tau; z)
&:=&\sqrt{\frac{y}{zy'}} \frac{1}{\sqrt{\alpha_1 \alpha_2}}
\frac{1}{n}
\sum_{t=1}^{\lfloor n\tau\rfloor} \prod_{i=1}^2 \big(1 - \frac{1}{\alpha_i}\big)^{t-s_i} \1(t \ge s_i), \qquad s_1, s_2 \in \Z, \\
\alpha_1&:=&\lambda' n y/z, \qquad  \alpha_2 \ := \ \lambda' n/y'. \nn
\end{eqnarray}
If
\begin{equation} \label{alphacond}
n, \alpha_1, \alpha_2, y, y' \to \infty \quad \text{so that} \quad  y/y' \to 1 \quad \text{and} \quad n = o(\alpha_i), \ i=1,2,
\end{equation}
then
\begin{eqnarray}
\widetilde h^{(\alpha_1,\alpha_2)}_{\alpha_1,\alpha_2}(s_1,s_2;\tau; z)
&:=&\sqrt{\alpha_1 \alpha_2} h_{\alpha_1,\alpha_2}(\lfloor \alpha_1 s_1\rfloor, \lfloor\alpha_2 s_2\rfloor; \tau; z) \nn \\
&=&\sqrt{\frac{y}{zy'}} \,\frac{1}{n}  \sum_{t=1}^{\lfloor n\tau\rfloor} \prod_{i=1}^2 \big(1 - \frac{1}{\alpha_i}\big)^{t-\lfloor\alpha_i s_i\rfloor}
\1(t \ge \lfloor \alpha_i s_i \rfloor) \nn \\
&\to&\frac{\tau}{\sqrt{z}}\prod_{i=1}^2 \e^{s_i} \1(s_i<0) =: h(s_1,s_2; \tau; z)\label{halphaconv}
\end{eqnarray}
point-wise for any $\tau >0$, $z > 0$, $s_i \in \R$, $s_i \ne 0$, $i=1,2 $ fixed.
Moreover, under the same conditions
\eqref{alphacond},
$\|\widetilde h^{(\alpha_1,\alpha_2)}_{\alpha_1,\alpha_2}(\cdot;\tau; z) - h(\cdot; \tau; z)\| \ \to \ 0 $, implying
the convergence $ Q_{12} (h_{\alpha_1,\alpha_2}(\cdot; \tau;z)) \to_{\rm d} I_{12} (h(\cdot; \tau; z)) =_{\rm d}
\tau Z_1 Z_2/(2\sqrt{z})$, $Z_i \sim N(0,1), i=1,2 $ by Proposition \ref{Qijconv}. Conditions
on $n, y, y', \lambda'$ in 
\eqref{alphacond} are obviously satisfied due to $y, y' = O(\lambda) = o(\lambda') $.
This proves \eqref{Llim} and \eqref{Vlim}, thereby completing the proof of
of \eqref{E1}.

\smallskip

\noi {\it Proof of \eqref{E2}.} For $T_{ij}(\theta)$ defined by \eqref{Tijdef} let us prove
\eqref{Tij}. Denote $N'_\lambda := (N \log \lambda)^{1/2\beta}$.
Similarly to \eqref{ineq:Tijbdd} we have that
\begin{eqnarray*} \label{Tbdd}
	|T_{ij}(\theta)|&\le&\frac{C}{N^{1/2} (\log \lambda)^{3/2}} \int_{(0, N'_\lambda]^3} \min \big\{ 1, \E z^2_{N,n} (\tau; x_1,x_2) \big\} (x_1 x_2 x_3)^{\beta-1} \1 (x_3 < x_1) \d x_1 \d x_2 \d x_3
\end{eqnarray*}
with $z_{N,n}(\tau; x_1,x_2)$ defined by \eqref{zndef31}. Whence using \eqref{Qnewbdd}
similarly as in the proof of case (i) we obtain
\begin{eqnarray*}
&|T_{ij}(\theta)|\ \le \
\frac{C}{N^{1/2} (\log \lambda)^{3/2}} \int_{(0, N'_\lambda]^2}
\min \Big\{1, \frac{1}{x_1 x_2} \min \big\{1, \frac{\lambda'}{x_1+x_2}\big\}\Big\} x_1^{2\beta -1} x_2 ^{\beta-1} \d x_1 \d x_2 \\
&\hskip-5cm = \frac{C}{N^{1/2} (\log \lambda)^{3/2}} \sum_{i=1}^3 T_{\lambda,i},
\end{eqnarray*}
where
\begin{eqnarray*}
T_{\lambda,1}&:=&\int_{(0, N'_\lambda]^2} \1( x_1 + x_2 < \lambda') \min \big\{1, \frac{1}{x_1 x_2} \big\} x_1^{2\beta -1} x_2 ^{\beta-1}
\d x_1 \d x_2, \\
T_{\lambda,2}&:=&\int_{(0, N'_\lambda]^2} \1(x_1 x_2 (x_1 + x_2) < \lambda', x_1 + x_2 > \lambda') x_1^{2\beta -1} x_2 ^{\beta-1}
\d x_1 \d x_2, \\
T_{\lambda,3}&:= &\lambda' \int_{(0, N'_\lambda]^2}
\1(x_1 x_2 (x_1 + x_2) > \lambda', x_1 + x_2 > \lambda' ) x_1^{2\beta -2} x_2 ^{\beta-2} \d x_1 \d x_2/ (x_1 + x_2).
\end{eqnarray*}
By changing variables $x_1, x_2$ as in \eqref{L1}-\eqref{Vint} we get
$T_{\lambda,1} \le C \int_0^{\lambda'} y^{\beta -1} \d y
\le C (\lambda')^\beta $. Also, similarly to the estimation of $J'_\lambda, J''_\lambda $, following \eqref{L2}
we obtain
$T_{\lambda,2} + T_{\lambda,3} \le C (\lambda')^\beta \int_{\lambda'}^{2N'_\lambda} y^{-1} \d y
\le C (\lambda')^\beta \log (N'_\lambda/\lambda') $. Hence, we conclude that
 \begin{eqnarray*}
|T_{ij}(\theta)|&\le&
\frac{C (\lambda')^\beta \log (N'_\lambda/\lambda') }{N^{1/2} (\log \lambda)^{3/2}} \ \le \
\frac{C \log n}{n^\beta \log \lambda},
\end{eqnarray*}
proving \eqref{Tij} with any $0< \delta < \beta $.
This proves \eqref{E2}. We omit the proof of \eqref{S2} which is completely similar to that in case (iii) and elsewhere.
This completes the proof of Theorem \ref{thmjoint} for $(t,s) = (0,1)$.

\smallskip

\noi {\it Proof of Theorem \ref{thmjoint} in the general case $(t,s) \in \Z^2, s \ge 1$.} Similarly
to \eqref{Sdecomp} we decompose $S^{t,s}_{N,n} (\tau)$ in \eqref{Sts} as
\begin{equation}\label{Stsdecomp}
S^{t,s}_{N,n}(\tau) = S^{t,s}_{N,n;q} (\tau) + S^{t,s; \dag}_{N,n;q} (\tau) + S^{t,s; \ddag}_{N,n;q} (\tau),
\end{equation}
where the main term
\begin{eqnarray}\label{Stsp}
S^{t,s}_{N,n;q} (\tau)&:=&\sum_{k=1}^{\tilde N_q} Y^{t,s}_{k,n;q} (\tau), \quad
Y^{t,s}_{k,n;q} (\tau) \ := \ \sum_{(k-1)q < i \le kq - s} \sum_{u=1}^{\lfloor n\tau\rfloor} X_i (u) X_{i+s}(u+t)
\end{eqnarray}
is a sum of independent $\tilde N_q=\lfloor N/q\rfloor$ blocks of size $q-s = q_{N,n} -s \to \infty$, and
\begin{eqnarray*}\label{Stsdag}
&&S^{t,s;\dag}_{N,n;q}(\tau)\ := \ \sum_{k=1}^{\tilde N_q} \sum_{kq - s < i \le kq} \sum_{u=1}^{\lfloor n\tau\rfloor} X_i (u) X_{i+s}(u+t), \quad
S^{t,s;\ddag}_{N,n;q} (\tau)\ := \ \sum_{q \tilde N_q < i \le N} \sum_{u=1}^{\lfloor n\tau\rfloor} X_{i} (u) X_{i+s}(u+t)
\end{eqnarray*}
are remainder terms. The proof of \eqref{E1}--\eqref{E2} for
$A^{-1}_{N,n}Y^{t,s}_{1,n;q}(\tau)
= \sum_{i=1}^{q-s} y^{t,s}_i(\tau), \ y^{t,s}_i(\tau) := A^{-1}_{N,n}\sum_{u=1}^{\lfloor n\tau\rfloor} X_i (u) X_{i+s}(u+t)$ is completely analogous since the distribution of $y^{t,s}_i(\tau) $ does not depend on $t$ and $s \neq 0$.

\subsection{Proof of Theorem \ref{thmclt}}

The proof uses the following result of \cite{pils2017}.

\begin{lemma}{(\cite{pils2017}, Lemma 7.1)}
	Let $\{ \xi_{ni}, 1 \le i \le N_n\}$, $n \ge 1$, be a triangular array of $m$-dependent r.v.s with zero mean and finite variance. Assume that: (L1) $\xi_{ni}$, $1 \le i \le N_n$, are identically distributed for any $n \ge 1$, (L2) $\xi_{n1} \to_{\rm d} \xi$, $\E \xi_{n1}^2 \to \E \xi^2 < \infty$ for some r.v.\ $\xi$ and (L3) $\operatorname{var} ( \sum_{i=1}^{N_n} \xi_{ni} ) \sim \sigma^2 N_n$, $\sigma^2 > 0$. Then $N_n^{-1/2} \sum_{i=1}^{N_n} \xi_{ni} \to_{\rm d} N(0,\sigma^2)$.
\end{lemma}

For notational simplicity, we consider only one-dimensional convergence at $\tau > 0$. Let $(Nn)^{-1/2} S_{N,n}^{t,s} (\tau) = N^{-1/2} \sum_{i=1}^N \xi_{ni}$, where $\xi_{ni} := n^{-1/2} \sum_{u=1}^{\lfloor n \tau \rfloor} X_i (u) X_{i+s} (u+t)$, $1 \le i \le N$ are $|s|$-dependent, identically distributed random variables with zero mean and finite variance.
Since $\xi_{ni}$, $1 \le i \le N$ are uncorrelated, it follows that
$\E ( \sum_{i=1}^N \xi_{ni} )^2 = N \E \xi_{n1}^2$, where $\xi_{n1} =_{\rm d} \xi_n := n^{-1/2} \sum_{u=1}^{\lfloor n\tau \rfloor} X_1(u) X_2 (u)$. Proposition \ref{propvar} implies $\E [ \xi_{n}^2 | a_1, a_2] \sim \tau A_{12}$, and so $\E \xi_n^2 \sim \tau \sigma^2$, where $\sigma^2 := \E A_{12} < \infty$.
It remains to show that $\xi_n \to_{\rm d} \sqrt{A_{12}} B(\tau)$, where $A_{12}$ is independent of $B(\tau)$. This follows from the martingale CLT similarly to \eqref{zlim}. By the lemma above, we conclude that $(Nn)^{-1/2} S^{t,s}_{Nn} (\tau) \to_{\rm d} \sigma B(\tau)$. Theorem \ref{thmclt} is proved. \hfill $\Box$

\smallskip

\section{Asymptotic distribution of temporal (iso-sectional) sample covariances}\label{sec:horiz}

The limit distribution of iso-sectional sample covariances $\widehat \gamma_{N,n}(t,0)$ in \eqref{covts}
and the corresponding partial sums process $S^{t,0}_{N,n}(\tau)$ of \eqref{Sts} is obtained similarly as in the
cross-sectional case, with certain differences which are discussed below. Since the conditional
expectation $\E [S^{t,0}_{N,n}(\tau)|a_1, \cdots, a_N] =: T^{t,0}_{N,n}(\tau) \ne 0$,
a natural decomposition is
\begin{equation}\label{St0decomp}
S^{t,0}_{N,n}(\tau)\ =\ \widetilde S^{t,0}_{N,n}(\tau)+ T^{t,0}_{N,n}(\tau),
\end{equation}
where  $\widetilde S^{t,0}_{N,n}(\tau) := S^{t,0}_{N,n}(\tau) - T^{t,0}_{N,n}(\tau)$
is the conditionally centered term
with $\E [\widetilde S^{t,0}_{N,n}(\tau)|a_1, \cdots, a_N] = 0$, and
\begin{equation}\label{Tt0}
T^{t,0}_{N,n}(\tau)\ :=\ \lfloor n\tau\rfloor\sum_{i=1}^N a^t_i/(1-a^2_i), \qquad t \ge 0,
\end{equation}
is proportional to a sum of i.i.d.\ r.v.s $a^t_i/(1-a^2_i), 1\le i \le N$ with regularly decaying tail distribution function
$$
\P \big( a^t/(1-a^2) > x) \sim \P (a > 1 - \frac{1}{2x}) \sim c_a x^{-\beta}, \qquad x \to \infty, \qquad
c_a := \psi(1)/2^\beta \beta,
$$
see condition \eqref{a:beta}. Accordingly,
the limit distribution of appropriately normalized and centered term $T^{t,0}_{N,n}(\tau)$ does not depend on $t$ and
can be found from the classical CLT and turns out
to be a $(\beta \wedge 2)$-stable line, under normalization $n N^{1/(\beta \wedge 2)}\, (\beta \ne 2). $
The other term, $\widetilde S^{t,0}_{N,n}(\tau)$, in \eqref{St0decomp}, is a sum of mutually independent partial sums processes
$Y^{t,0}_{i,n}(\tau) := \sum_{u=1}^{\lfloor n\tau \rfloor} (X_i(u) X_i(u+t) - \E [X_i(u) X_i(u+t)|a_i]), 1\le i \le N$ with
conditional variance
\begin{eqnarray*}
&{\rm var}[Y^{t,0}_{i,n} (1)|a_i] \sim n A^{t,0}_{ii}, \quad
n \to \infty, \qquad \text{where} \quad A^{t,0}_{ii} := \frac{1+a_i^2}{1-a_i^2} \Big( \frac{1 + a_i^{2|t|}}{(1-a_i^2)^2} + \frac{a_i^{2|t|} (2|t| + \operatorname{cum}_4)}{1-a_i^4} \Big).
\end{eqnarray*}
The proof of the last fact follows similarly to that of \eqref{Aijdef} and is omitted.
As $a_i \uparrow 1 $, $A^{t,0}_{ii} \sim 1/(2(1-a_i)^3)$ and the limit
distribution of $\widetilde S^{t,0}_{N,n}(\tau)$ can be shown to exhibit a trichotomy on the interval $0< \beta < 3 $
depending on the limit $ \lambda^\ast_\infty $ in \eqref{mucond2}. It turns out that
for $\beta > 2 $ the asymptotically Gaussian term $T^{t,0}_{N,n}(\tau)$ dominates $\widetilde S^{t,0}_{N,n}(\tau)$ in all
cases of $ \lambda^\ast_\infty $, while in the interval $0< \beta < 2 $
$T^{t,0}_{N,n}(\tau)$ 
and $\widetilde S^{t,0}_{N,n}(\tau)$ have the same convergence rate. 
Somewhat surprisingly, the limit distribution of $S^{t,0}_{N,n}(\tau)$
is a $\beta$-stable line in both cases $  \lambda^\ast_\infty = \infty$ and $  \lambda^\ast_\infty = 0$ with
different scale parameters of the random slope coefficient of this line.

Rigorous description of the above limit results is given in the following Theorems \ref{thmjoint2} and \ref{thmclt2}. The proofs
of these theorems are similar and actually simpler than the corresponding
Theorems \ref{thmjoint} and \ref{thmclt} dealing with non-horizontal sample covariances, due to
the fact that $S^{t,0}_{N,n}(\tau)$ is a sum of row-independent summands contrary to $S^{t,s}_{N,n}(\tau), s \ne 0$.
Because of this,
we omit some details of the proof of Theorems \ref{thmjoint2} and \ref{thmclt2}. We also omit the more delicate
cases $\beta =1 $ and $\beta = 2 $ where the limit results may require a change of normalization or additional
centering.

\begin{theorem} \label{thmjoint2} Let the mixing distribution satisfy condition \eqref{a:beta} with $0< \beta <2$, $\beta\ne 1$. Let
$N, n \to \infty $ so that
\begin{equation}\label{mucond2}
\lambda^\ast_{N,n} := \frac{N^{1/\beta}}{n} \ \to \lambda^\ast_\infty \ \in \ [0, \infty].
\end{equation}
In addition, assume $\E \vep^4(0) < \infty $.
Then the following statements (i)--(iii) hold for $S^{t,0}_{N,n}(\tau), t \in \Z$ in \eqref{Sts} depending on $\lambda^\ast_\infty$ in \eqref{mucond2}.

\medskip

\noi (i) Let $\lambda^\ast_\infty = \infty$. Then
\begin{eqnarray} \label{muinf2}
&&
n^{-1} N^{-1/\beta} 
\big(S^{t,0}_{N,n}(\tau) - \E S^{t,0}_{N,n}(\tau) \1(1< \beta < 2)\big)
\ \to_{\rm fdd} \ \tau V^\ast_\beta,
\end{eqnarray}
where
$V^\ast_{\beta}$ is a completely asymmetric $\beta$-stable r.v. with characteristic function in \eqref{Vstar} below.

\medskip

\noi (ii) Let $\lambda_\infty^\ast = 0$. Then
\begin{equation}\label{muzero2}
n^{-1} N^{-1/\beta} \big(S^{t,0}_{N,n}(\tau) - \E S^{t,0}_{N,n}(\tau) \1(1< \beta < 2)\big) \ \to_{\rm fdd} \
\tau V^+_\beta,
\end{equation}
where $V^+_{\beta}$ is a completely asymmetric $\beta$-stable r.v. with characteristic function in \eqref{Vplus} below.

\medskip

\noi (iii) Let $0<\lambda_\infty^* < \infty$. Then
\begin{eqnarray} \label{interlim2}
&&
n^{-1} N^{-1/\beta} 
\big(S^{t,0}_{N,n}(\tau) - \E S^{t,0}_{N,n}(\tau)\1(1< \beta < 2)\big) \ \to_{\rm fdd}\ \lambda_\infty^*
{\cal Z}^\ast_{\beta}(\tau/\lambda_\infty^*), 
\end{eqnarray}
where ${\cal Z}^\ast_{\beta} $ is the `diagonal intermediate' process in \eqref{interZ01}.

\end{theorem}

\begin{remark} The r.v.s $V^\ast_\beta $ and $V^+_\beta $
	in \eqref{muinf2} and \eqref{muzero2}
	have respective stochastic integral representations
	\begin{eqnarray*}
		V^\ast_\beta &=&\int_{\R_+ \times C(\R)} \big\{ \int_{-\infty}^0 \e^{xs} \d B(s) \big\}^2  \d ({\cal M}^\ast_\beta -  \E {\cal M}^\ast_\beta \1(1< \beta < 2)),\\
		V^+_\beta &=&\int_{\R_+ \times C(\R)} (2x)^{-1} \d ({\cal M}^\ast_\beta -  \E {\cal M}^\ast_\beta \1(1< \beta < 2))
	\end{eqnarray*}
	w.r.t. Poisson random measure ${\cal M}^\ast_\beta $
	in  \eqref{mubeta1}. Note $\int_{-\infty}^0 \e^{xs} \d B(s) =_{\rm law} Z/\sqrt{2x}, Z \sim N(0,1)$.
	The fact that both
	$ V^\ast_\beta $ and $V^+_\beta $ 
	have { completely asymmetric} $\beta$-stable distribution follows from
	their ch.f.s:
	{
		\begin{eqnarray}
		\E \e^{\i \theta V^\ast_\beta}
		&=&\exp \big\{\psi(1) \int_0^\infty \E \big(\e^{\i \theta Z^2/{(2x)}} - 1 - \i (\theta Z^2/{(2x)}) \1(1< \beta < 2)\big) x^{\beta -1} \d x \big\}\nn\\
		&=& \exp \big\{ -c^*_\beta |\theta|^\beta (1 - \i \operatorname{sign}(\theta) \tan (\pi \beta/2) ) \big\}, \label{Vstar}\\
		\E \e^{\i \theta V^+_\beta}
		&=& \exp \big\{\psi(1) \int_0^\infty \big(\e^{\i \theta/{(2x)}} - 1 - \i (\theta/{(2x)}) \1(1< \beta < 2)\big) x^{\beta -1} \d x \big\}\nn\\
		&=& \exp\big\{ - c^+_\beta |\theta|^\beta (1- \i \operatorname{sign}(\theta) \tan (\pi \beta/2) ) \big\}, \quad \theta \in \R,\label{Vplus}
		\end{eqnarray}
		where
		\begin{eqnarray}
		c^+_\beta := \frac{\psi (1)\Gamma(2-\beta)\cos (\pi \beta/2)}{2^\beta \beta (1-\beta)}, \quad c^\ast_\beta := c^+_\beta \E  |Z|^{2\beta} \label{castbeta}
		\end{eqnarray}
		with}
	$\E  |Z|^{2\beta} = 2^\beta \Gamma(\beta + 1/2)/\sqrt{\pi} \ne 1 $ unless
	$\beta = 1$, implying that  $V^\ast_\beta$ and  $V^+_\beta $ have different distributions.

\end{remark}

\begin{theorem} \label{thmclt2}
Let the mixing distribution satisfy condition \eqref{a:beta} with $\beta > 2 $. In addition, assume $\E \vep^4(0) < \infty $.
Then for any
$t \in \Z$,
as $N, n \to \infty $ in arbitrary way,
\begin{equation}
n^{-1} N^{-1/2} \big(S^{t,0}_{N,n} (\tau)- \E S^{t,0}_{N,n} (\tau) \big)\ \to_{\rm fdd}\ \tau \sigma^*_t Z, 
\end{equation}
where $Z \sim N(0,1)$ and
$(\sigma^*_t)^2 
:=  {\rm var} (a^{|t|}/(1-a^2))$.

\end{theorem}

\begin{remark}
If $\beta < 1$, then $\gamma(t,0)$ is undefined for any $t \in
\Z$. Using the convention $\gamma(t,0) \1(1< \beta < 2) := 0$ if $ \beta < 1$ and $\gamma(t,0)$ if $\beta > 1$.
\end{remark}

\begin{corollary} \label{cor2}

\noi (i) Let the conditions of Theorem \ref{thmjoint2} (i) be satisfied. Then for any $t \in \Z$
\begin{eqnarray*}
&N^{1 - 1/\beta} (\widehat \gamma_{N,n} (t,0) - \gamma(t,0) \1(1< \beta < 2)) \ \to_{\rm d} \ V^\ast_\beta.
\end{eqnarray*}

\smallskip

\noi (ii) Let the conditions of Theorem \ref{thmjoint2} (ii) be satisfied. Then for any $t \in \Z$
\begin{eqnarray*}
&N^{1 - 1/\beta} (\widehat \gamma_{N,n} (t,0) - \gamma(t,0) \1(1< \beta < 2)) \ \to_{\rm d} \ V^+_{\beta}.
\end{eqnarray*}

\noi (iii) Let the conditions of Theorem \ref{thmjoint2} (iii) be satisfied. Then for any $t \in \Z$
\begin{eqnarray*}
&N^{1 - 1/\beta} (\widehat \gamma_{N,n} (t,0) - \gamma(t,0) \1(1< \beta < 2)) \ \to_{\rm d} \ \lambda^*_\infty {\cal Z}^*_\beta (1/\lambda^*_\infty).
\end{eqnarray*}

\smallskip

\noi (iv) Let the conditions of Theorem \ref{thmclt2} be satisfied. Then for any $t \in \Z$
\begin{eqnarray*}
&N^{1/2} (\widehat \gamma_{N,n} (t,0) - \gamma(t,0)) \ \to_{\rm d} \ \sigma^*_t Z, \qquad Z \sim N(0,1).
\end{eqnarray*}

\end{corollary}

\noi {\it Proof of Theorem \ref{thmjoint2}.} Let $t\ge 0$ and 
\begin{eqnarray}
y^{t,0}(\tau)&:=&\frac{1}{n N^{1/\beta}}\sum_{u=1}^{\lfloor n\tau\rfloor} (X(u)X(u+t)- \E X (u) X(u+t) \1(1<\beta <2)).
\end{eqnarray}
It suffices to prove that
\begin{eqnarray}\label{Phiqconv1}
&\Phi^{t,0}_{N,n}({\theta}) \ \rightarrow \ \Phi^\ast ({\theta}), \quad \text{as} \ N, \, n \to \infty, \, \lambda^\ast_{N,n} \to
\lambda^\ast_\infty, \quad \forall \theta \in \R,
\end{eqnarray}
where, using $\E y^{t,0}(\tau) \1(1< \beta < 2) = 0$,
\begin{eqnarray} \label{Phip2}
&\Phi^{t,0}_{N,n}(\theta)
:= N \E \big[\e^{\i \theta y^{t,0}(\tau) } -1 - \i \theta y^{t,0}(\tau) \1( 1< \beta < 2) \big], \qquad
\Phi^\ast(\theta) := \log \E \e^{\i \theta {\cal S}^\ast_{\beta} (\tau)},
\end{eqnarray}
and ${\cal S}^\ast_{\beta} (\tau)$ denotes the limit process in \eqref{muinf2}--\eqref{interlim2}.
Similarly to \eqref{Phirep},
\begin{eqnarray}
\Phi^{t,0}_{N,n}(\theta)
&=&\psi(1)\int_{(0,1/N^{1/\beta}]}
\E \big[\e^{\i \theta z^{t,0}_{N,n}(\tau;x)} -1 - \i \theta z^{t,0}_{N,n}(\tau;x) \1(1< \beta < 2) \big] x^{\beta-1} \d x,
\label{Phirep2}
\end{eqnarray}
where $z^{t,0}_{N,n}(\tau;x) := y^{t,0}(\tau)|_{a = 1 - x/N^{1/\beta}}$.
Next we decompose $y^{t,0}(\tau) = y^{\ast}(\tau) + y^{+}(\tau)$, where
\begin{eqnarray*}
y^{\ast}(\tau)&:=&\frac{1}{n N^{1/\beta}}\sum_{u=1}^{\lfloor n\tau\rfloor} (X(u)X(u+t)- \E [X (u) X(u+t)|a]), \\
y^{+}(\tau)&:=&\frac{\lfloor n\tau\rfloor}{n N^{1/\beta}}(\E [X (0) X(t)|a] - \E [X (0) X(t) \1(1< \beta < 2)]) \ =\
\frac{\lfloor n\tau\rfloor}{n N^{1/\beta}}\Big(\frac{a^t}{1-a^2} - \E \Big[\frac{a^t\1(1< \beta < 2)}{1-a^2}\Big]\Big).
 \end{eqnarray*}
Accordingly, we decompose $ z^{t,0}_{N,n}(\tau;x)= z^{\ast}_{N,n}(\tau;x) + z^{+}_{N,n}(\tau;x)$, where
\begin{eqnarray} \label{zndef2}
z^{\ast}_{N,n}(\tau;x)
&:=&\frac{1}{n N^{1/\beta}} \sum_{s_1,s_2 \in \Z} \overline{\vep(s_1) \vep(s_2)}
\sum_{u=1}^{\lfloor n\tau\rfloor} \big(1 - \mbox{$\frac{x}{N^{1/\beta}}$}\big)^{2u+t-s_1-s_2}\1(u \ge s_1,u+t \ge s_2), \\
z^+_{N,n}(\tau;x)
&:=&\frac{\lfloor n\tau\rfloor}{n N^{1/\beta}}
\Big(\frac{(1 - x N^{-1/\beta})^t} {1 - (1- x N^{-1/\beta})^2 } -
\E \Big[\frac{a^t\1(1< \beta < 2)}{1-a^2}\Big]\Big), \nn
\end{eqnarray}
where $\overline{\vep(s_1) \vep(s_2)}
:= \vep(s_1) \vep(s_2)- \E\vep(s_1) \vep(s_2). $

\medskip

\noi {\it Proof of \eqref{Phiqconv1}, case $0< \lambda^\ast_\infty < \infty$. } We have
\begin{eqnarray} \label{Phip3}
&\Phi^\ast(\theta)= \psi(1) 
\int_0^\infty
\E \big[\e^{\i \theta \lambda^*_\infty z^\ast(\tau/\lambda^*_\infty; x)} -1 - \i \theta
\lambda^*_\infty z^\ast(\tau/\lambda^*_\infty; x) \1(1 < \beta < 2)\big] x^{\beta -1} \d x,
\end{eqnarray}
where the last expectation is taken w.r.t.\ the Wiener measure $P_B $.
Similarly as in the proof of \eqref{E1} we
prove 
the point-wise convergence of the integrands in \eqref{Phirep2} and \eqref{Phip3}: for any $x > 0$
\begin{eqnarray} \label{zast}
\Lambda^{t,0}_{N,n}(\theta; x)&:=&\E \big[\e^{\i \theta z^{t,0}_{N,n}(\tau;x) } -1 -
\i \theta z^{t,0}_{N,n}(\tau; x) \1(1< \beta < 2)
\big]\\
&\to&\E \big[\e^{\i \theta \lambda^*_\infty z^{\ast}(\tau/\lambda^*_\infty;x) } -1 - \i \theta \lambda^*_\infty
z^{\ast}(\tau/\lambda^*_\infty;x) \1(1< \beta < 2) \big].\nn
\end{eqnarray}
The proof of
\eqref{zast} using Proposition \ref{Qijconv} is very similar to that of \eqref{Qnconv} and we omit the details.
Using \eqref{zast} and the dominated convergence theorem we can prove the convergence of integrals, or \eqref{Phiqconv1}. The application of the dominated convergence theorem is guaranteed by the dominating bound
\begin{eqnarray}\label{Lambdabdd}
&|\Lambda^{t,0}_{N,n}(\theta;x)| \le C (1 \wedge (1/x))\{1(0 < \beta < 1)+ (1/x)\1(1< \beta < 2)\},
\end{eqnarray}
which is a consequence of $|z^+_{N,n}(\tau;x)| \le C/x,
\E (z^{\ast}_{N,n}(\tau;x))^2
\le C x^{-2}$, see
\eqref{var1}. Particularly, for $0< \beta < 1$ we get $|\Lambda_{N,n}^{t,0}(\theta;x)| \le 2$ and $|\Lambda_{N,n}^{t,0}(\theta;x)| \le
\E (|z^{\ast}_{N,n}(\tau;x)| + |z^+_{N,n}(\tau;x)|) \le C ( \sqrt{\E |z^{\ast}_{N,n}(\tau;x)|^2} + (1/x) ) \le C/x $, hence \eqref{Lambdabdd} follows. For $1< \beta < 2$ \eqref{Lambdabdd} follows similarly.
This proves \eqref{Phiqconv1} for $0< \lambda^\ast_\infty < \infty$.

\medskip

\noi {\it Proof of \eqref{Phiqconv1}, case $\lambda^\ast_\infty = 0$. } In this case
\begin{eqnarray*} \label{Phip4}
&\Phi^\ast(\theta)= \psi(1) \int_{\R_+}
\big[\e^{\i \theta (\tau/(2 x))} -1 - \i\theta (\tau/(2x)) \1(1 < \beta < 2)\big] x^{\beta -1} \d x,
\end{eqnarray*}
see \eqref{Vplus}.
From \eqref{var1} we have $\E (z^{\ast}_{N,n}(\tau;x))^2
\le C x^{-2} \min \{1, \lambda^\ast_{N,n}/x\} \to 0 $ and hence
\begin{eqnarray} \label{zast2}
\Lambda^{t,0}_{N,n}(\theta; x)
&\to&\e^{\i \theta \tau/(2 x)} -1 - \i \theta
(\tau/(2 x)) \1(1< \beta < 2)\nn
\end{eqnarray}
for any $x >0$ similarly as in \eqref{zast}. Finally, the use of the dominating bound in
\eqref{Lambdabdd} which is also valid in this case completes the proof of
 \eqref{Phiqconv1} for $\lambda^\ast_\infty = 0$.

\medskip

\noi {\it Proof of \eqref{Phiqconv1}, case $\lambda^\ast_\infty = \infty$. } In this case,
\begin{eqnarray} \label{Phip5}
&\Phi^\ast(\theta)= \psi(1) \int_{\R_+}
\E \big[\e^{\i \theta (\tau Z^2/(2x))} -1 - \i\theta (\tau Z^2/(2x)) \1(1 < \beta < 2)\big] x^{\beta -1} \d x,
\end{eqnarray}
see \eqref{Vstar}. Write $z^{\ast}_{N,n}(\tau;x)$ in
\eqref{zndef2} as quadratic form: $z^{\ast}_{N,n}(\tau;x) = Q_{11} (h(\tau; x; \cdot)) $ in
\eqref{Qij} and apply Proposition \ref{Qijconv} with $\alpha_1 = \alpha_2 \equiv \alpha := N^{1/\beta} $.
Note $\widetilde h^{(\alpha, \alpha)} (\tau; x; s_1, s_2) =
n^{-1} \sum_{u=1}^{\lfloor n\tau\rfloor} (1 - x/N^{1/\beta})^{u-\lfloor N^{1/\beta} s_1\rfloor}(1 - x/N^{1/\beta})^{t+u- \lfloor N^{1/\beta} s_2\rfloor}
\1(u \ge \lfloor N^{1/\beta} s_1\rfloor, u+t \ge \lfloor N^{1/\beta} s_2\rfloor) \to g(s_1,s_2) := \tau \e^{x(s_1+s_2)} \1(s_1 \vee s_2 \le 0)
$ point-wise a.e.\ in $(s_1, s_2) \in \R^2$ and also in $L^2(\R^2)$. Then conclude
$z^{\ast}_{N,n}(\tau;x) \to_{\rm d} I_{11} (g) =_{\rm d} \int_{\R^2} g(s_1,s_2) \d B(s_1) \d B(s_2)
=_{\rm d} \tau \big\{\big(\int_{-\infty}^0 \e^{sx} \d B(s) \big)^2 - \E \big(\int_{-\infty}^0 \e^{sx} \d B(s) \big)^2 \big\}
=_{\rm d} \tau (Z^2 - 1)/(2x) $ for any $x>0$,
 where $Z \sim N(0,1)$. On the other hand,
$z^+_{N,n}(\tau;x) \to \tau/2x $ and therefore
\begin{eqnarray*} \label{zast3}
\Lambda^{t,0}_{N,n}(\theta; x)
&\to&\E \big[\e^{\i \theta \tau Z^2/(2 x)} -1 - \i \theta
(\tau Z^2/(2 x)) \1(1< \beta < 2)\big]\nn
\end{eqnarray*}
for any $x >0$,
proving 
the point-wise convergence of the integrands in \eqref{Phirep2} and \eqref{Phip5}. The remaining details
are similar as in the previous cases and omitted. This ends the proof of Theorem \ref{thmjoint2}.
\hfill $\Box$

\medskip

\noi {\it Proof of Theorem \ref{thmclt2}.} Consider the decomposition
in \eqref{St0decomp}, where $n^{-1} T^{t,0}_{N,n}(\tau) = (\lfloor n\tau\rfloor/n) \sum_{i=1}^N a^t_i/(1-a^2_i)$
is a sum of i.i.d.\ r.v.s with finite variance $(\sigma_t^*)^2 =  {\rm var} (a^{|t|}/(1-a^2))$ and therefore
\begin{equation*}
n^{-1} N^{-1/2} \big(T^{t,0}_{N,n} (\tau)- \E T^{t,0}_{N,n} (\tau) \big)\ \to_{\rm fdd}\ \tau \sigma^*_t Z
\end{equation*}
holds by the classical CLT as $N, n \to \infty $ in arbitrary way and where $Z \sim N(0,1)$. Hence, the statement of the theorem follows from
$\widetilde S^{t,0}_{N,n} (1) = o_{\rm p}( n N^{1/2})$.
By Proposition \ref{propvar} \eqref{var1} we have that
$ {\rm var} (\widetilde S^{t,0}_{N,n}(1)) = N \E {\rm var} [\sum_{u=1}^n X(u)X(u+t) |a ]
\le  C N n^2 \E \big[ (1-a)^{-2} \min \{1, (n (1-a))^{-1} \} \big] $, where
the last expectation vanishes as $n \to \infty $, due to $\E (1-a)^{-2} < \infty $.
Theorem \ref{thmclt2} is proved. \hfill $\Box$

\begin{figure}[h!]
	\begin{tabular}{cr}
		\includegraphics[width=0.5\textwidth]{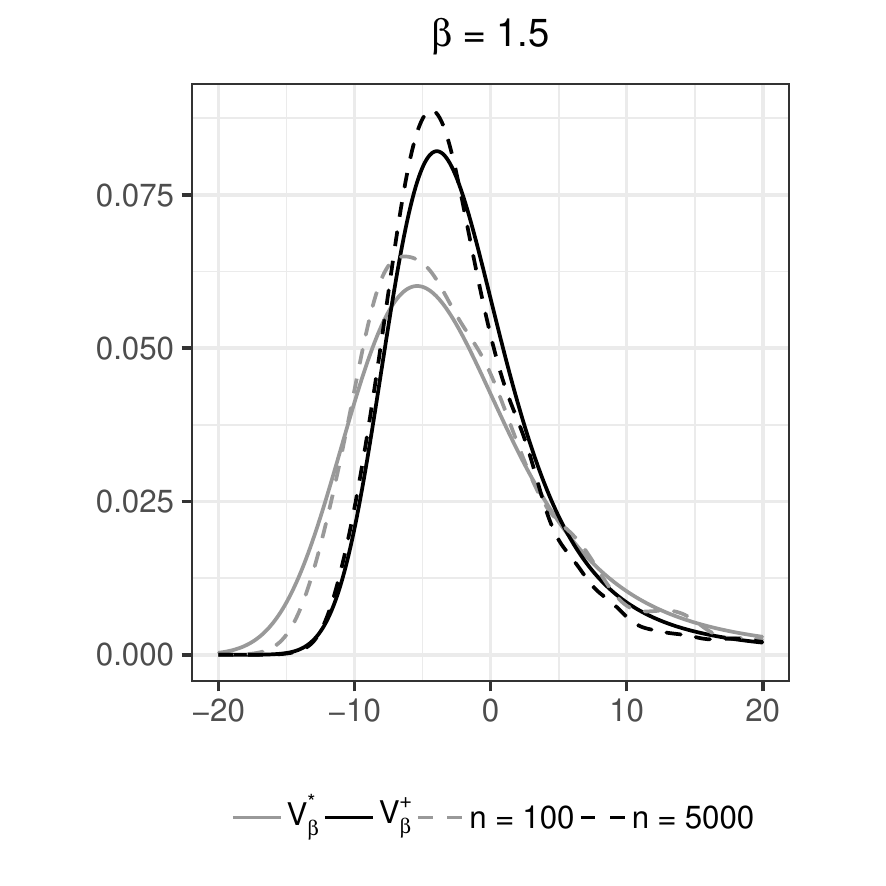}	
		&\includegraphics[width=0.5\textwidth]{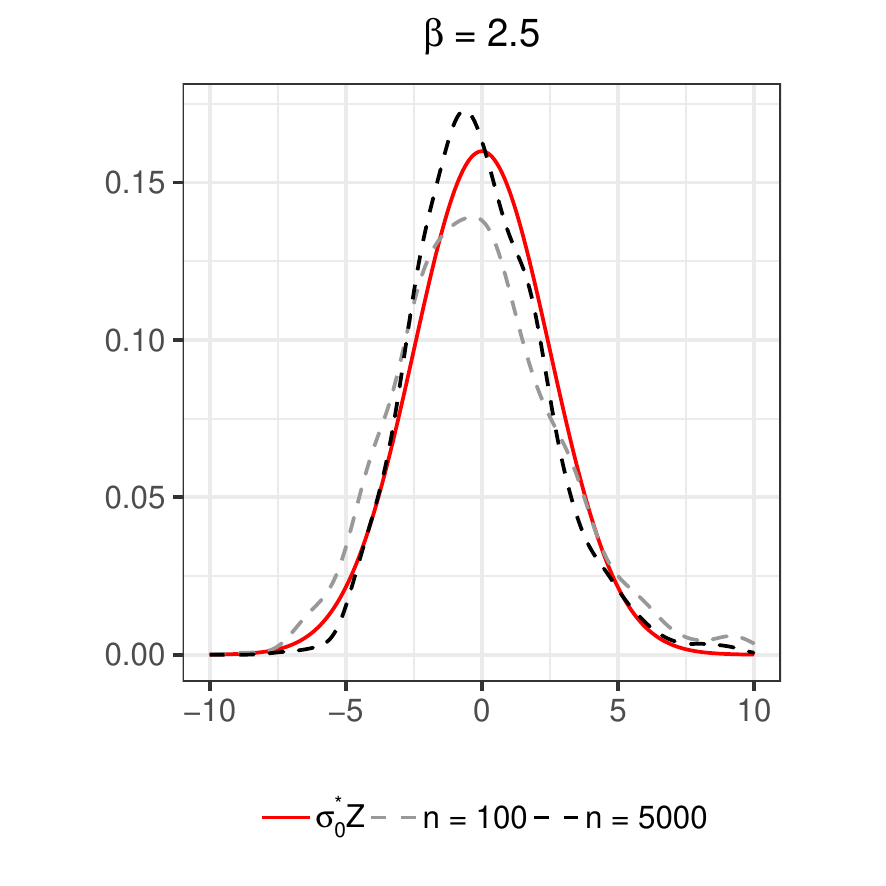}
	\end{tabular}
	\caption{Density of the limiting random variables in cases [left] (i),(ii), [right] (iv) of Corollary \ref{cor2} for $t=0$ and their kernel density estimates constructed from a random sample of size 1000 from $\widehat \gamma_{N,n} (0,0) 
		$ in  \eqref{cov00} with $N=5000$, $a^2 \sim \operatorname{Beta}(2,\beta)$, $\vep(0) \sim N(0,1)$.
	}
	\label{fig}
\end{figure}

{
To illustrate our results, we use $a^2 \sim \operatorname{Beta}(\alpha,\beta)$, $\alpha, \beta > 0$, as in \cite{gran1980}.
Then condition \eqref{a:beta} holds with the same $\beta$ and we can explicitly compute parameters of the limit distributions in cases (i), (ii), (iv) of Corollary~\ref{cor2}.} Figure \ref{fig} shows the density of the corresponding limiting random variables for $\alpha=2$, $\beta=1.5$, $2.5$ and $t=0$. We also plot the kernel density estimates constructed using
1000 RCAR(1) panels with $N=5000$, $n=100$, $5000$, $\vep(0) \sim N(0,1)$.
More specifically, we use a random sample of $N^{1/\beta} (\widehat \gamma_{N,n} (0,0) -\gamma(0,0))$ if $\beta=1.5$ and $N^{1/2} (\widehat \gamma_{N,n} (0,0) -\gamma(0,0))$ if $\beta = 2.5$.
On the l.h.s.\ we can see that the empirical distribution of $\widehat \gamma_{N,n} (0,0)$ is different for $n=100$, $5000$, whereas on the r.h.s.\ both kernel density estimates are quite close to the limiting normal density.

{ In the finite variance case $\beta > 1 $,
Corollary \ref{cor2} can be used for statistical inference about the covariance  $\gamma (t,0) = \gamma (t)$ in \eqref{gammat}, provided
parameters of the limit distributions are consistently estimated. Denote by
\begin{eqnarray}\label{FFdf}
F^*_{\beta,\psi}(x) := \P(V^\ast_\beta \le x), \quad
F^+_{\beta,\psi}(x) := \P(V^+_\beta \le x),  \qquad x \in \R,
\end{eqnarray}
the c.d.f.s of the above stable r.v.s, which are uniquely determined by $\beta$, $\psi (1) \equiv \psi $ in \eqref{a:beta}, see \eqref{Vstar}--\eqref{castbeta}.
The same is true for the (marginal) distribution ${\cal Z}^*_\beta (\tau)$ of the `diagonal intermediate' process in
\eqref{interZ01}.
In Corollary \ref{cor3}  we suppose the existence of estimators
\begin{eqnarray}\label{betapsi}
\hat \beta_{N,n}&=&\beta + o_{\rm p}(1/{\log N}), \qquad \hat  \psi_{N,n} = \psi + o_{\rm p}(1), \\
\hat \sigma^2_{N,n,t}&=&(\sigma^*_{t})^2 + o_{\rm p}(1), \label{sigmat}
\end{eqnarray}
which is discussed in Remark \ref{rem4} below.
Corollary \ref{cor3} omits the `intermediate' case  $\lambda^*_\infty \in (0,\infty)$, partly because  in this case
the limit distribution
is less tractable and depends on  $\lambda^*_\infty $  which is difficult to assess in a finite sample.

\begin{corollary} \label{cor3}
	(i) Let the conditions of Theorem \ref{thmjoint2} (i) be satisfied, $1< \beta < 2 $, and $\hat \beta_{N,n}, \hat \psi_{N,n}$
	be estimators as in \eqref{betapsi}. Then for any $t \in \Z$
	\begin{eqnarray} \label{sup1}
	&\sup_{x \in \R} \big|\P \big(N^{1 - 1/\hat \beta_{N,n}} (\widehat \gamma_{N,n} (t,0)  - \gamma(t)) \le x \big)
	- F^*_{\hat \beta_{N,n}, \hat \psi_{N,n}}(x)\big|
	\  = \ o_{\rm p}(1).
	\end{eqnarray}

	\noi (ii) Let the conditions of Theorem \ref{thmjoint2} (ii) be satisfied, $1< \beta < 2$,
	and $\hat \beta_{N,n}, \hat \psi_{N,n}$
	be estimators as in \eqref{betapsi}. Then for any $t \in \Z$
	\begin{eqnarray} \label{sup2}
	&\sup_{x \in \R} \big|\P \big(N^{1 - 1/\hat \beta_{N,n}} (\widehat \gamma_{N,n} (t,0)  - \gamma(t)) \le x \big)
	- F^+_{\hat \beta_{N,n}, \hat \psi_{N,n}}(x)\big|
	\  = \ o_{\rm p}(1).
	\end{eqnarray}

	\noi (iii) Let the conditions of Theorem \ref{thmclt2} be satisfied, $\beta > 2$,
	and $\hat \sigma^2_{N,n,t}$
	be an estimator as in \eqref{sigmat}. Then for any $t \in \Z$
	\begin{eqnarray} \label{sup3}
	&\sup_{x \in \R} \big|\P \big(\big(\frac{N}{ \hat \sigma^2_{N,n,t}}\big)^{1/2}
	(\widehat \gamma_{N,n} (t,0)  - \gamma(t)) \le x \big)
	- \P(Z\le x)\big|
	\  = \ o_{\rm p}(1), \qquad Z \sim N(0,1).
	\end{eqnarray}
	
\end{corollary}

\noi {\it Proof.} Consider \eqref{sup1}.
Write $N^{1 - 1/\hat \beta_{N,n}} (\widehat \gamma_{N,n} (t,0)  - \gamma(t))
= N^{1 - 1/\beta} (\widehat \gamma_{N,n} (t,0)  - \gamma(t)) + \xi_{N,n}$, where
$\xi_{N,n} := (N^{{(1/\beta)-(1/\hat \beta_{N,n})}} -1) N^{1 - 1/\beta} (\widehat \gamma_{N,n} (t,0)  - \gamma(t)) = o_{\rm p}(1) $ due to
\eqref{betapsi} and Corollary \ref{cor2}(i). Therefore,
$\sup_{x \in \R} |\P(N^{1 - 1/\hat \beta_{N,n}} (\widehat \gamma_{N,n} (t,0)  - \gamma(t)) \le x )
- F^*_{\beta, \psi}(x) | \to 0$. Relation $\sup_{x \in \R} | F^*_{\beta, \psi}(x)
- F^*_{\hat \beta_{N,n}, \hat \psi_{N,n}}(x) | = o_{\rm p}(1)$ follows from \eqref{betapsi} and continuity of continuity of the c.d.f. $F^*_{\beta,\psi}$
in $\beta, \psi$. {This proves \eqref{sup1}. The proof of \eqref{sup2}, \eqref{sup3} is analogous.} \hfill $\Box$

\begin{remark} \label{rem3} Using Corollary \ref{cor3} we can construct asymptotic confidence intervals
	for $\gamma (t)$, as follows.
	For $\alpha \in (0,1)$ denote by $q_{\beta,\psi}(\alpha)$ the $\alpha$-quantile of the c.d.f. $F^*_{\beta,\psi}$ in \eqref{FFdf}. Then, since $\alpha = \allowbreak F^*_{\hat \beta_{N,n}, \hat \psi_{N,n}}( q_{\hat \beta_{N,n}, \hat \psi_{N,n}} (\alpha) )$ a.s., $\P(
	N^{1 - 1/\hat \beta_{N,n}} (\widehat \gamma_{N,n} (t,0)  - \gamma(t)) \le q_{\hat \beta_{N,n},\hat \psi_{N,n}} (\alpha))
	- \alpha = o_{\rm p} (1)$
	follows from \eqref{sup1}; moreover since the above quantity is non-random, we get that $|\P (N^{1 - 1/\hat \beta_{N,n}} (\widehat \gamma_{N,n} (t,0)  - \gamma(t)) \le q_{\hat \beta_{N,n},\hat \psi_{N,n}} (\alpha) )
	- \alpha | = o(1)$,
	implying that
	\begin{eqnarray*}
		&\Big[\widehat \gamma_{N,n} (t,0)  - N^{(1/\hat \beta_{N,n})-1} q_{\hat \beta_{N,n},\hat \psi_{N,n}} (1-\alpha/2), \widehat \gamma_{N,n} (t,0)  - N^{(1/\hat \beta_{N,n})-1} q_{\hat \beta_{N,n},\hat \psi_{N,n}} (\alpha/2) \Big]
	\end{eqnarray*}
is the asymptotic
confidence interval for $\gamma(t)$, for any confidence level $\alpha \in (0,1)$. Analogous confidence intervals for
$\gamma(t)$ can be defined in the case \eqref{sup2}; in the case   \eqref{sup3} they follow in a standard way.

\end{remark}

\begin{remark} \label{rem4} Estimation of the tail parameter $\beta $ in the RCAR(1) panel model was studied
	in \cite{lei2019}. Particularly, \cite{lei2019} developed a modified version $\hat \beta_{N,n}$
	of the Goldie--Smith estimator in \cite{gol1987} and proved its asymptotic normality, under additional
	(rather stringent) conditions on the mutual increase rate of $N $ and $n$. A similar estimator   $\hat \psi_{N,n}$  can be
	defined following \cite{gol1987}. We expect that
	these estimators satisfy the consistency as in  \eqref{betapsi} under much weaker assumptions on  $N, n$.
	Finally, for $t \ge 0$ the estimator $\hat \sigma^2_{N,n,t}$  in \eqref{sigmat} can be defined (see the proof in Appendix)
as
		\begin{eqnarray} \label{hsigmat}
		\hat \sigma^2_{N,n,t}
		&:=&\frac{1}{N} \sum_{i=1}^N \Big(\frac{1}{n} \sum_{k=1}^{n-t} X_i(k)X_i(k+t) \Big)^2 - \Big(\frac{1}{Nn} \sum_{i=1}^N \sum_{k=1}^{n-t} X_i(k)X_i(k+t)\Big)^2.
		\end{eqnarray}
	\end{remark}

\begin{remark}
	In general, in the RCAR(1) model the  autoregressive coefficient $a$ can take a value from $(-1,1)$.
 In the latter case if the distribution of $a$ is sufficiently dense at $-1$, the (unconditional) autocovariance function of the RCAR(1) process oscillates when decaying slowly, which is usually referred to as seasonal long memory.
The restriction $a \in [0,1) $ in the present paper (as well as in \cite{pils2017}, \cite{lei2019} and some other papers)
is basically due to technical reasons. We expect that, under assumption  \eqref{a:beta}, %
most of our results hold in the general case $a \in (-1,1) $ provided
the concentration of the mixing distribution near $-1$ is not too strong, e.g., if   $\E (1+a)^{-\beta'} < \infty$ for some $\beta'>\beta$.

	
\end{remark}

}
\appendix

\section{Appendix}


\noi {\it Proof of Proposition \ref{propinter}}. (i) The existence of ${\cal Z}_\beta $ follows from
\begin{equation}\label{Lc}
J_\beta := \int_{{\cal L}_1^c } |z(\tau; x_1,x_2)|^2 \d \mu_\beta < \infty
\end{equation}
and $\mu_\beta ({\cal L}_1) < \infty. $
We have $\mu_\beta ({\cal L}_1) = \psi(1)^2 \int_{\R^2_+} \1(x_1 x_2 (x_1 + x_2) < 1) (x_1 x_2)^{\beta -1} \d x_1 \d x_2
\le C \int_0^\infty x_1^{\beta -1} \d x_1 $ $\int_0^{x_1} \1 ( x_2 < 1/x_1^2) x_2^{\beta -1} \d x_2 $ $ = C \big(\int_0^1 x_1^{\beta -1} \d x_1
\int_0^{x_1} x_2^{\beta -1} \d x_2 + \int_1^\infty x_1^{\beta -1} \d x_1 \int_0^{1/x_1^2} x_2^{\beta -1} \d x_2 \big)
\le C \big(\int_0^1 x_1^{2\beta -1} \d x_1 + \int_1^\infty x_1^{-\beta -1} \d x_1 \big) < \infty
$ since $\beta > 0$.

Consider \eqref{Lc}. Then
$$
J_\beta = C \int_{\R^2_+} \1(x_1 x_2 (x_1 + x_2) > 1) \E |z(\tau; x_1,x_2)|^2 (x_1 x_2)^{\beta -1}
\d x_1 \d x_2,
$$
where
\begin{align}\label{z2}
\E |z(\tau; x_1,x_2)|^2
&=\int_{(0,\tau]^2} \prod_{i=1}^2 \E [{\cal Y}_i(u_1;x_i){\cal Y}_i(u_2;x_i)] \d u_1 \d u_2\nn\\
&= \frac{1}{4 x_1 x_2} \int_{(0,\tau]^2} \e^{-(x_1+x_2)|u_1-u_2|} \d u_1 \d u_2
\le \frac{C\tau^2}{x_1x_2} \Big( 1 \wedge \frac{1}{\tau(x_1+x_2)} \Big).
\end{align}
Hence,
\begin{align*}
J_\beta &\le C \int_{\R^2_+} \1 (x_1x_2(x_1+x_2) > 1 ) (x_1+x_2)^{-1} (x_1x_2)^{\beta-2} \d x_1 \d x_2\\
&\le C \int_{\R^2_+} \1 (x_2 > x_1, \, x_1x_2^2 > 1 ) x_1^{\beta-2} x_2^{\beta-3} \d x_1 \d x_2\\
&=C \Big( \int_0^1 x_1^{\beta-2} \d x_1 \int_{x_1^{-1/2}}^\infty x_2^{\beta-3} \d x_2 + \int_1^\infty x_1^{\beta-2} \d x_1 \int_{x_1}^\infty x_2^{\beta-3} \d x_2 \Big) < \infty
\end{align*}
if $0 < \beta < 3/2$.
The remaining facts in (i) are easy and we omit the details.

\smallskip

\noi (ii) 
Similarly as in (\cite{pils2014}, proof of Proposition~3.1(ii))
it suffices to show for any $0< p < 2\beta$ that
\begin{equation}
\infty > J_{p,\beta}(\tau):=
\begin{cases}
\int_{\R^2_+} \E |z(\tau; x_1,x_2)|^p (x_1 x_2)^{\beta -1} \d x_1 \d x_2, &0< p \le 2, \\
\int_{\R^2_+} \E \big[|z(\tau; x_1,x_2)|^p \vee |z(\tau; x_1,x_2)|^2 \big]
 (x_1 x_2)^{\beta -1} \d x_1 \d x_2, &p> 2.
\end{cases}
 \label{Jp}
\end{equation}
Let first $0< p \le 2 $. Using $\E |z(\tau;x_1,x_2) |^p \le ( \E |z(\tau; x_1,x_2)|^2 )^{p/2}$ and \eqref{z2}, we obtain
\begin{align}
J_{p,\beta} (\tau) \le C \int_{\R^2_+} \Big( \int_{(0,\tau]^2} \e^{-(x_1+x_2)|u_1-u_2|} \d u_1 \d u_2 \Big)^{p/2} (x_1x_2)^{\beta-1-p/2} \d x_1 \d x_2=: C \tau^{2(p-\beta)} I_{p,\beta}\label{Iptau},
\end{align}
where
\begin{align}
I_{p,\beta} &\le \int_{\R^2_+} \Big( 1 \wedge \frac{1}{x_1+x_2} \Big)^{p/2} (x_1 x_2)^{\beta-1-p/2} \d x_1 \d x_2\nn\\
&\le C \int_0^\infty \int_0^{x_1} \Big( 1 \wedge \frac{1}{x_1} \Big)^{p/2} (x_1x_2)^{\beta-1-p/2} \d x_1 \d x_2\nn \\
&= C \int_0^\infty \Big( 1 \wedge \frac{1}{x_1} \Big)^{p/2} x_1^{2\beta - p - 1} \d x_1 < \infty\label{Ipb}
\end{align}
if $p/2 < \beta < 3p/4$, 
thus proving \eqref{Jp} for $0< p \le 2 $.

Next for $ 2 < p < 3 $ we need the inequality for double It\^o-Wiener integrals: for any $p\ge 2, g \in L^2(\R^2)$
\begin{eqnarray}\label{I2g}
\textstyle \E \big|\int_{\R^2} g(s_1,s_2) \d B_1 (s_1) \d B_2 (s_2) \big|^p
\le C \big( \E \big|\int_{\R^2} g(s_1,s_2) \d B_1 (s_1) \d B_2 (s_2) \big|^2\big)^{p/2}
= C \big(\int_{\R^2} |g(s_1,s_2)|^2 \d s_1 \d s_2\big)^{p/2}.
\end{eqnarray}
Indeed, by using Gaussianity and independence of $B_1, B_2$ and Minkowski inequality for $I_2(g) := \int_{\R^2} g(s_1,s_2) \d B_1 (s_1) \d B_2 (s_2)$
we obtain
\begin{eqnarray*}
\big(\E | I_2(g)|^p \big)^{2/p}
&=&\big(\E_{B_1} \E_{B_2} \big[ |I_2(g)|^p \big| B_1 \big]\big)^{2/p}
\le C\big( \E_{B_1} (\E_{B_2} \big[ |I_2(g)|^2 \big| B_1 \big])^{p/2} \big)^{2/p} \\
&\le&C\E_{B_2} \big\{\E_{B_1} \big[ |I_2(g)|^p \big| B_2 \big]\}^{2/p} \le
C\E_{B_2} \big\{ \big(\E_{B_1} \big[ |I_2(g)|^2 \big| B_2 \big] \big)^{p/2} \big\}^{2/p} \\
&=&C\E_{B_2} \E_{B_1} \big[ |I_2(g)|^2 \big| B_2 \big] = C \E |I_2(g)|^2.
\end{eqnarray*}
Using inequality \eqref{I2g} and \eqref{Iptau}, \eqref{Ipb} we obtain
\begin{align*}
J_{p,\beta}(\tau) &\le C\Big(\int_{\R^2_+}\E |z(\tau; x_1,x_2)|^p (x_1 x_2)^{\beta -1} \d x_1 \d x_2 
+ \int_{\R^2_+} \E |z(\tau; x_1,x_2)|^2 (x_1 x_2)^{\beta -1} \d x_1 \d x_2 \Big)\\
&\le C(I_{p,\beta} (\tau)+ I_{2,\beta}(\tau)) < \infty
\end{align*}
if $p/2 < \beta < 3p/4$, thus proving \eqref{Jp} and part (ii).

\smallskip

\noi (iii) Follows from stationarity of increments of ${\cal Z}_\beta $ (part (i)) and
$J_{2,\beta} (\tau) = \sigma^2_\infty \tau^{2(2-\beta)} $, where
according to \eqref{z2},
\begin{eqnarray}
\sigma^2_\infty&=&\int_{\R^2_+} \E z^2(1;x_1,x_2) \d \mu_\beta \nn \\
&=&(\psi(1)/2)^2\int_{(0,1]^2} \d u_1 \d u_2 \big(\int_0^\infty \e^{-x |u_1-u_2|} x^{\beta -2} \d x \big)^2 \nn \\
&=&(\psi(1)/2)^2 \Gamma(\beta -1)^2 \int_{(0,1]^2} |u_1-u_2|^{2(1-\beta)} \d u_1 \d u_2 \
= \  (\psi(1)/2)^2 \Gamma(\beta -1)^2/((2-\beta)(3-2\beta)). \nn \label{sigmainf}
\end{eqnarray}

\smallskip

\noi (iv) Follows from stationarity of increments,
$\E |{\cal Z}_\beta(\tau)|^p \le CJ_{p,\beta} (\tau)$, $1 < p \le 2$, where $J_{p,\beta} (\tau)$ is the same as in \eqref{Jp}, and Kolmogorov's criterion;
c.f (\cite{pils2014}, proof of Proposition 3.1(iv)).

\smallskip

\noi (v) The proofs are very similar to those of Theorem \ref{thmjoint} (i), (ii), hence we omit some details. For notational simplicity, we only prove one-dimensional convergence at $\tau > 0$.

\smallskip

\noi {\it Proof of \eqref{ass:local1}.} As $b \to 0$, consider
$$
\Phi_b (\theta) := \log \E \exp \{ \i \theta b^{\beta-2} {\cal Z}_\beta (b \tau) \} = \psi(1)^2 \int_{\R_+^2} \E \Psi (\theta b^{\beta-2} z(b\tau; x_1,x_2) ) (x_1 x_2)^{\beta-1} \d x_1 \d x_2,
$$
where $\Psi(z) := \e^{\i z} - 1 - \i z$, $z \in \R$. Since
$b^{-2} z(b\tau;x_1,x_2) =_{\rm d} z(\tau; bx_1, bx_2)$,
rewrite
$$
\Phi_b (\theta) = \psi(1)^2 b^{-2\beta} \int_{\R^2_+} \E \Psi (\theta b^\beta z(\tau;x_1,x_2)) (x_1x_2)^{\beta-1} \d x_1 \d x_2,
$$
where $b^{-2\beta} \Psi (\theta b^\beta z(\tau;x_1,x_2)) \to - (\theta^2/2) z^2 (\tau; x_1, x_2)$ a.s. Note $|b^{-2\beta} \Psi (\theta b^\beta z(\tau; x_1,x_2))| \le (\theta^2/2) z^2 (\tau; x_1,x_2)$, where the dominating function satisfies \eqref{z2} and \eqref{Ez5}. Hence, by the dominated convergence theorem,
$$
\Phi_b (\theta) \to - (\theta^2/2) \psi(1)^2 \int_{\R_+^2} \E z^2 (\tau;x_1,x_2) (x_1 x_2)^{\beta-1} \d x_1 \d x_2 = \log \E \{ \i \theta \sigma_\infty B_{2-\beta} (\tau) \},
$$
which finishes the proof.

\smallskip

\noi {\it Proof of \eqref{ass:local0}} follows that of Thm.\ \ref{thmjoint}(i), case $0 < \beta < 1$. As $b \to 0$, consider
\begin{align*}
\Phi_b (\theta) := \log \E \e^{ \i \theta b^{-1} (\log b^{-1})^{-1/2\beta} {\cal Z}_\beta ( b \tau ) } = \frac{\psi(1)^2}{\log b^{-1}} \int_{\R^2_+} \E [ \e^{ \i \theta z_b (\tau;x_1,x_2) } - 1 ] (x_1x_2)^{\beta-1} \d x_1 \d x_2,
\end{align*}
where
$$
z_b (\tau;x_1,x_2) := b^{-1} (\log b^{-1})^{-1/(2\beta)} z \big( b \tau ; ( \log b^{-1} )^{-1/(2\beta)} x_1, ( \log b^{-1} )^{-1/(2\beta)} x_2 \big)
$$
satisfies
\begin{equation}\label{ineq:varzb}
\E |z_b (\tau;x_1,x_2)|^2 \le \frac{C}{x_1 x_2} \Big(1 \wedge \frac{b^{-1} (\log b^{-1})^{1/(2\beta)}}{x_1 + x_2} \Big),
\end{equation}
see \eqref{z2}.
Split
\begin{align*}
\Phi_b (\theta) &= \frac{\psi(1)^2}{\log b^{-1}} \int_{\R^2_+} \big( \1 (1 < x_1+x_2 < b^{-1} ) + \1 (x_1+x_2> b^{-1} ) + \1 (x_1+x_2 < 1) \big)\\
&\quad \times \E [ \e^{ \i \theta z_b (\tau;x_1,x_2) } - 1 ] (x_1x_2)^{\beta-1} \d x_1 \d x_2 =: \sum_{i=1}^3 L_i.
\end{align*}
Using \eqref{ineq:varzb}, we can show that $L_i$, $i=2,3$ are remainders.
By change of variables: $y = x_1 + x_2, x_1 = yw $ and then $w = z/y^2$,
we rewrite the main term
\begin{equation}\label{L1b}
L_1 = \frac{1}{\log b^{-1}} \int_1^{b^{-1}} V_b (\theta;y) \frac{\d y}{y}, \qquad
V_b (\theta;y) := 2 \psi(1)^2 \int_0^{y^2/2} \Lambda_b (z;y) z^{\beta-1} \big( 1-\frac{z}{y^2} \big)^{\beta-1} \d z
\end{equation}
with
$\Lambda_b (z;y) := \E [ \exp \{ \i \theta z_b (\tau; \frac{z}{y}, y (1-\frac{z}{y^2} ) ) \} - 1 ]$, which
satisfies $|\Lambda_b (z;y)| \le C (1 \wedge \frac{{1}}{z})$ for all $0 < \frac{z}{y^2} < \frac{1}{2}$, $0 < y < b^{-1}$. Here the dominating bound is a consequence of \eqref{ineq:varzb}.
Then
\begin{equation}\label{Vlimb}
L_1 \to \log \E \e^{ \i \theta \tau V_{2\beta} } = 2 \psi(1)^2 \int_0^{\infty}\Lambda (z) z^{\beta -1} \d z,
\end{equation}
where $\Lambda (z) := \E [\e^{\i \theta \tau Z_1 Z_2/(2\sqrt{z})} -1]$ with $Z_i \sim N(0,1)$, $i=1,2$ being independent r.v.s, follows from
\begin{eqnarray}\label{Llimb}
\lim_{y \to \infty, y = O(b^{-1})} \Lambda_b (z; y) = \Lambda (z), \qquad \forall z > 0,
\end{eqnarray}
for more details we refer the reader to the proof of Thm.\ \ref{thmjoint} (i) case $0< \beta < 1$. More precisely, \eqref{Llimb} says that for every $\epsilon >0 $ there exists a small $\delta>0$ such that for all $0< b < \delta$, if $\delta^{-1} < y < b^{-1}$, then $|\Lambda_b (z; y) - \Lambda(z)| < \epsilon$.
To show \eqref{Llimb}, note $z_b ( \tau;\frac{z}{y}, y (1-\frac{z}{y^2} ) ) = I_{12} (h_b (\cdot;\tau;z) )$ is a double It{\^o}-Wiener stochastic integral w.r.t.\ independent standard Brownian motions $\{B_i(s), s\in \R \}$, $i=1,2$ for
\begin{align*}
&h_b (s_1,s_2;\tau;z) := (\log b^{-1})^{-1/(2\beta)} \int_0^{\tau} \prod_{i=1}^2 \e^{- \frac{1}{\alpha_i} (bu-s_i)} \1 (s_i < bu) \d u, \quad s_1,s_2 \in \R,\\
&\alpha_1 := (\log b^{-1})^{1/(2\beta)} y/z, \quad
\alpha_2 := (\log b^{-1})^{1/(2\beta)}/ y', \quad y' := y \big( 1-\frac{z}{y^2} \big).
\end{align*}
We have that
$z_b (\tau;\frac{z}{y}, y(1-\frac{z}{y^2})) =_{\rm d} I_{12} ( \widetilde h_b ( \cdot ;\tau;z) )$, where
\begin{align*}
\widetilde h_b (s_1,s_2;\tau;z) &:= \sqrt{\alpha_1 \alpha_2} h_b (\alpha_1 s_1, \alpha_2 s_2; \tau; z)\\
&= \sqrt{\frac{y}{z y'}} \int_0^{\tau} \prod_{i=1}^2 \e^{- \frac{1}{\alpha_i} (bu - \alpha_i s_i)} \1 (\alpha_i s_i < bu) \d u, \quad s_1,s_2 \in \R.
\end{align*}
If $b \to 0$, $y, y' \to \infty$ so that $y/y' \to 1$ and $b/\alpha_i \to 0$, $i=1,2$,
then
$\|\widetilde h_b (\cdot;\tau; z) - h(\cdot; \tau; z)\| \to 0$ with
\begin{eqnarray}
h(s_1,s_2; \tau; z) := \frac{\tau}{\sqrt{z}}\prod_{i=1}^2 \e^{s_i} \1(s_i<0), \quad s_1, s_2 \in \R,
\end{eqnarray}
implies the convergence $z_b (\tau;\frac{z}{y}, y(1-\frac{z}{y^2})) \to_{\rm d} I_{12} (h(\cdot; \tau; z)) =_{\rm d}
\tau Z_1 Z_2/2\sqrt{z}$. Conditions
on $b, y, y'$ are obviously satisfied due to $y, y' = O(b^{-1}) = o( b^{-1} (\log b^{-1})^{1/(2\beta)} )$.
This proves \eqref{Llimb} and \eqref{Vlimb}, thereby completing the proof of
of \eqref{ass:local0}.

\smallskip

\noi {\it Proof of \eqref{ass:global}} follows that of Theorem \ref{thmjoint} (ii). We will prove that as $b \to \infty$,
\begin{align}
\log \E \e^{\i \theta b^{-1/2} {\cal Z}_\beta (b \tau) } &= \psi (1)^2 \int_{\R^2_+} \E \big[ \exp \big\{ \i \theta b^{-1/2} z (b \tau;x_1,x_2) \big\} - 1 \big] (x_1 x_2)^{\beta-1} \d x_1 \d x_2\label{Phirepb}\\
&\to \psi(1)^2 \int_{\R_+^2} \Big[ \exp \Big\{ - \frac{\theta^2 \tau}{4 x_1 x_2 (x_1 + x_2)} \Big\} - 1 \Big] (x_1 x_2)^{\beta-1} \d x_1 \d x_2 = \log \E \e^{\i \theta {\cal A}^{1/2} B(\tau)}.\nn
\end{align}
By \eqref{z2}, we have that $\E [ \exp \{ \i \theta b^{-1/2} z (b \tau;x_1,x_2) \} - 1] \le C \min \{1, (x_1 x_2 (x_1 + x_2))^{-1}\}$. In view of \eqref{Aineq2}, the dominated convergence theorem applies if the integrands on the r.h.s.\ of \eqref{Phirepb} converge pointwise, i.e.\ for every $(x_1,x_2) \in \R^2_+$,
\begin{equation}\label{zlimb}
b^{-1/2} z (b\tau; x_1, x_2) \to_{\rm d} \frac{B(\tau)}{\sqrt{2 x_1 x_2 (x_1 + x_2)}}.
\end{equation}
To simplify notation, let $\tau =1 $ and all $b \in \N$. Define
$$
z^+_b (x_1,x_2)
:= \int_0^{b} \int_0^{b} f (s_1,s_2) \d B_1(s_1) \d B_2 (s_2), \quad f (s_1,s_2) := b^{-1/2} \int_0^{b} \prod_{i=1}^2 \e^{-x_i(u-s_i)} \1 (u > s_i) \d u,
$$
and
$z^-_b (x_1,x_2) := b^{-1/2} z (b;x_1,x_2) - z^+_b (x_1,x_2)$. Since
$\E (z^-_b (x_1,x_2))^2 = O (b^{-1})$ implies $z^-_b (x_1,x_2) = o_{\rm p} (1)$, we only need to prove that
\begin{equation} \label{zzlimb}
z^+_b (x_1,x_2) \to_{\rm d} N \Big(0, \frac{1}{2 x_1 x_2 (x_1 + x_2)} \Big) \quad \text{as } b\to\infty.
\end{equation}
Write $z^+_b (x_1,x_2)=\sum_{k=1}^b Z_k$ as a sum of a sum of a zero-mean square-integrable martingale difference array
\begin{align*}
Z_k &:= \int_{k-1}^{k} \int_{0}^{k-1} f (s_1,s_2) \d B_1 (s_1) \d B_2 (s_2) + \int_{0}^{k-1} \int_{k-1}^{k} f (s_1,s_2) \d B_1 (s_1) \d B_2 (s_2)\\
&\quad+ \int_{k-1}^k \int_{k-1}^{k} f (s_1,s_2) \d B_1 (s_1) \d B_2 (s_2)
\end{align*}
w.r.t.\ the filtration ${\cal F}_k$ generated by $\{ B_i (s), 0 \le s \le k, \, i =1,2 \}$, $k=0, \dots, b$.
By the martingale CLT in Hall and Heyde \cite{hall1980}, \eqref{zzlimb} then follows from
\begin{equation}\label{a:CLTb}
\sum_{k=1}^b \E [Z^2_{k} | {\cal F}_{k-1}] \to_{\rm p} \frac{1}{2 x_1 x_2 (x_1+x_2)} \quad \text{and} \quad
\sum_{k=1}^b \E [Z^2_{k} \1 (|Z_{k}|> \epsilon )] \to 0 \quad \text{for any } \epsilon > 0.
\end{equation}
Since $\sum_{k=1}^b \E Z^2_k = \int_0^b \int_0^b f^2 (s_1,s_2) \d s_1 \d s_2 = \E (z^+_b (x_1,x_2))^2 \to (2 x_1 x_2 (x_1+x_2) )^{-1}$, consider $R_b := \sum_{k=1}^b ( \E [Z^2_k | {\cal F}_{k-1} ] - \E Z^2_k )$, where
\begin{align*}
\E [ Z^2_k | {\cal F}_{k-1} ] &= \int_{k-1}^k \Big( \int_0^{k-1} f(s_1,s_2) \d B_2 (s_2) \Big)^2 \d s_1 + \int_{k-1}^k \Big( \int_0^{k-1} f(s_1,s_2) \d B_1 (s_1) \Big)^2 \d s_2\\
&\quad + \int_{k-1}^k \int_{k-1}^k f^2 (s_1,s_2) \d s_1 \d s_2.
\end{align*}
By rewriting
$R_b =_{\rm d} \sum_{i=1}^2 \int_0^b \int_0^b c_i (s_1,s_2) \d B_i (s_1) \d B_i (s_2)$ with
$c_1 (s_1,s_2) = \int_{\lceil s_1 \vee s_2 \rceil}^b f(s_1,s) f(s_2,s) \d s$,
$c_2 (s_1,s_2) = \int_{\lceil s_1 \vee s_2 \rceil}^b f(s,s_1) f(s,s_2) \d s$
and using the elementary bound:
\begin{equation}\label{ineq:f}
f(s_1,s_2) \le C b^{-1/2} \big( \e^{-x_1 (s_2-s_1)} \1 (s_1 < s_2) + \e^{-x_2(s_1-s_2)} \1 (s_1 \ge s_2) \big), \quad 0 \le s_1, s_2 \le b,
\end{equation}
we obtain
$\E |R_b|^2 = \sum_{i=1}^2 \int_0^b \int_0^b c^2_i(s_1,s_2) \d s_1 \d s_2 = O(b^{-1}) = o(1)$, which proves $R_b = o_{\rm p} (1)$ and completes the proof of the first relation in \eqref{a:CLTb}.
Finally, using \eqref{I2g}, \eqref{ineq:f}, we obtain $\sum_{k=1}^b \E |Z_k|^4= O(b^{-1}) = o (1)$, which implies the second relation in \eqref{a:CLTb} and completes the proof of \eqref{zzlimb}.

\noi Proposition \ref{propinter} is proved. \hfill $\Box$

\medskip

\noi {\it Proof of Proposition \ref{propinter1}.} (i) Split ${\cal Z_\beta^\ast} (\tau) = \widetilde {\cal Z}_\beta^\ast (\tau) + \tau V^+_\beta$ with
\begin{align*}
	\widetilde {\cal Z}_\beta^\ast (\tau) &:= \int_{\R_+ \times C(\R)} \Big( z^\ast (\tau;x) - \frac{\tau}{2x} \Big) \d ( {\cal M}_\beta^\ast - \E {\cal M}_\beta^\ast \1(1<\beta <2) ),\\
	V_\beta^+ &:= \int_{\R_+ \times C(\R)} \frac{1}{2x} \d ( {\cal M}_\beta^\ast - \E {\cal M}_\beta^\ast \1(1<\beta <2) ),
\end{align*}
where ${\cal M}^*_\beta$ is a Poisson random measure on $\R_+ \times C(\R)$ with mean $\mu_\beta^* = \E {\cal M}^*_\beta$ given in \eqref{mubeta1}.
The existence of $V_\beta^+$ follows from $\int_0^\infty \min\{1, x^{-1}\} x^{\beta-1} \d x < \infty$ if $\beta \in (0,1)$ and $\int_0^\infty \min\{x^{-1}, x^{-2}\} x^{\beta-1} \d x < \infty$ if $\beta \in (1,2)$. The process $\widetilde {\cal Z}_\beta^\ast$ is well-defined if
\begin{equation}\label{Jstar}
J^\ast_{p,\beta} (\tau) := \int_{\R_+ \times C(\R)} |z^\ast (\tau;x) - \tau/2x |^p \d \mu_\beta^\ast = C \int_0^\infty \E |z^\ast (\tau;x) - \tau/2x |^p x^{\beta-1} \d x < \infty,
\end{equation}
where $0 < p \le 1$ for $\beta \in (0,1)$ and $1 \le p \le 2$ for $\beta \in (1,2)$. We have $\E |z^\ast(\tau;x) - \tau/2x|^p \le (\operatorname{var} (z^\ast (\tau;x)) )^{p/2}$, where
\begin{align}
	\operatorname{var} ( z^\ast (\tau;x) ) &= \int_{(0,\tau]^2} {\rm cov} ({\cal Y}^2(u_1;x), {\cal Y}^2(u_2;x)) \d u_1 \d u_2\nn\\
	&=2 \int_{(0,\tau]^2} \int_{\R^2} \d s_1 \d s_2 \e^{-2x(u_1+ u_2 - s_1 - s_2)} \1(s_1 \vee s_2 < u_1 \wedge u_2)\nn\\
	&= \frac{1}{2x^2} \int_{(0,\tau]^2} \e^{-2x|u_1-u_2|} \d u_1 \d u_2 = \frac{1}{8x^4} (2x\tau - 1 +e^{-2x\tau}) \le C \frac{\tau^2}{x^2} \Big( 1 \wedge \frac{1}{x\tau} \Big),\label{ineq:varzstar}
\end{align}
hence, $J^\ast_{p,\beta} (\tau) \le C \tau^{2p - \beta} < \infty$ for $p < \beta < 3p/2$. This completes the proof of part (i).

\smallskip

\noi (ii) $\E |V^+_\beta|^p < \infty$ for $0 < p<\beta$, since $V^+_\beta$ is a $\beta$-stable random variable. Similarly to \eqref{Jp}, $\E |\widetilde {\cal Z}_\beta^\ast (\tau)|^p < \infty$ follows from $J^*_{p,\beta} (\tau) < \infty$ in \eqref{Jstar}, where $p$ is sufficiently close to $\beta$ and such that $0 < p < \beta < 3p/2$.
This proves part (ii).

\smallskip	

\noi (iii) Follows from part (ii) by Kolmogorov's criterion, similarly as in the proof of Proposition \ref{propinter}.

\smallskip

\noi (iv) For notational simplicity, we only prove one-dimensional convergence at $\tau > 0$. We have
$\log \E \exp \{ \i \theta b^{-1} {\cal Z}_\beta (b \tau) \} = \psi (1)\int_{\R_+} \Lambda_b (x) x^{\beta-1} \d x$, where
$$
\Lambda_b (x) := \E \big[ \exp \big\{ \i \theta b^{-1} z^\ast (b \tau; x) \big\} - 1 - \i \theta b^{-1} z^\ast (b \tau; x) \1 (1 < \beta < 2) \big].
$$
Substituting $\E |z^\ast (b\tau; x)| \le (\E |z^\ast (b\tau; x)|^2)^{1/2}$ and
$\E|z^\ast (b\tau; x)|^2 = \operatorname{var} ( z^\ast (b\tau; x) ) + (b\tau/2x)^2 \le C (b/x)^2
$ by \eqref{ineq:varzstar} into
\begin{align*}
|\Lambda_b (x)| \le C
\begin{cases}
\min \big\{ 1, b^{-1} \E |z^\ast (b \tau;x)| \big\}, &0 < \beta < 1,\\
\min \big\{ b^{-1} \E |z^\ast (b \tau;x)|, b^{-2} \E |z^\ast (b \tau;x)|^2 \big\}, &1 < \beta < 2,
\end{cases}
\end{align*}
we obtain the bounds: $|\Lambda_b (x)| \le C \min \{ 1, x^{-1} \}$ if $0 < \beta < 1$, and $|\Lambda_b (x)| \le C \min \{ x^{-1}, x^{-2} \}$ if $1 < \beta < 2$.
The result then follows from the dominated convergence theorem once we show that for all $x \in \R_+$,
\begin{align}\label{lim:Lambdab}
\Lambda_b (x) \to
\begin{cases}
\exp \{\i \theta \tau/(2x) \} - 1 - (\i \theta \tau/(2x)) \1 (1 < \beta < 2) &\text{as } b \to \infty,\\
\E [ \exp \{\i \theta Z^2 \tau /(2x) \} - 1 - (\i \theta Z^2 \tau/(2x)) \1 (1 < \beta < 2) ] &\text{as } b \to 0,
\end{cases}
\end{align}
where $Z \sim N(0,1)$.
Using \eqref{ineq:varzstar}, we get $\E | b^{-1} z^\ast (b\tau; x) - (\tau/2x) |^2 = b^{-2} \operatorname{var} ( z^\ast (b \tau; x) ) \le C b^{-1}= o (1)$ as $b \to \infty$, which implies the first convergence in \eqref{lim:Lambdab}.
To prove the second convergence in \eqref{lim:Lambdab}, note $Z /\sqrt{2x} =_{\rm d} {\cal Y} (0; x)$. It suffices to show that as $b \to 0$,
\begin{align*}
\E |b^{-1} z^\ast (b\tau; x) - \tau {\cal Y}^2 (0; x) | = \E \Big| \int_0^\tau ( {\cal Y}^2 (b u; x) - {\cal Y}^2 (0; x) ) \d u \Big| \le \int_0^\tau \E | {\cal Y}^2 (b u; x) - {\cal Y}^2 (0; x) | \d u = o(1).
\end{align*}
Factorizing the difference of squares and applying the Cauchy-Schwarz inequality, this follows from
$$
\E | {\cal Y} (bu;x) - {\cal Y} (0;x) |^2 = \int_0^{bu} \e^{-2xs} \d s + \frac{1}{2x} (e^{- x b u} - 1)^2 \le C b u.
$$
Prop. \ref{propinter1} is proved.
\hfill $\Box$

\medskip

\noi {\it Calculation of the constant $\sigma_0$ in Proposition \ref{propinter} (v).}
We have
\begin{eqnarray*}
\sigma_0 \cdot \frac{2^{2\beta/3}}{\psi (1)^2} &=& \int_{\R_+^2} \big(1 - \exp \{ -(u_1+u_2)^{-1} (u_1 u_2)^{-1} \} \big) (u_1 u_2)^{\beta -1} \d u_1 \d u_2\\
&\underset{u_2 = u_1 v_2}{=} & \int_{\R_+^2} \big(1 - \exp \{ -u_1^{-3} (1+ v_2)^{-1} v_2^{-1} \} \big) u_1^{2\beta-1} v_2^{\beta -1} \d u_1 \d v_2\\
& \underset{u_1 = v_1^{-1/3}}{=}& \frac{1}{3} \int_{\R_+^2} \big(1 - \exp \{ -v_1 (1+ v_2)^{-1} v_2^{-1} \} \big) v_1^{-2\beta/3-1} v_2^{\beta -1} \d v_1 \d v_2\\
&= & \frac{1}{3} \int_{\R_+^2} \Big( \int_0^{1/((1+ v_2) v_2)} \e^{-v_1 t} \d t \Big)  v_1^{-2\beta/3} v_2^{\beta -1} \d v_1 \d v_2\\
&= & \frac{\Gamma (1-\frac{2\beta}{3}) }{3} \int_0^\infty v_2^{\beta -1} \d v_2
\int_0^{1/((1+ v_2) v_2)} t^{2\beta/3-1} \d t \\
&=& \frac{\Gamma (1-\frac{2\beta}{3}) }{2\beta} \int_0^\infty (1+ v_2)^{-2\beta/3} v_2^{\beta/3-1} \d v_2\\
& \underset{v_2 = s^{-1}-1}{=}& \frac{\Gamma (1-\frac{2\beta}{3}) }{2\beta} \int_0^1 s^{2\beta/3} (s^{-1}-1)^{\beta/3-1} s^{-2} \d s\\
&=& \frac{\Gamma (1-\frac{2\beta}{3}) \operatorname{B} (\frac{\beta}{3}, \frac{\beta}{3})}{2\beta} .
\end{eqnarray*}

{
\noindent {\it Proof of \eqref{sigmat}.} By Corollary~4.1~(iv), $\frac{1}{Nn} \sum_{i=1}^N \sum_{k=1}^{n-t} X_i(k) X_i(k+t)
 \to_{\rm p}  \gamma (t) = \E \frac{a^t}{1-a^2} $. Hence, relation \eqref{sigmat} for \eqref{hsigmat}
 follows from
\begin{equation} \label{X1}
\frac{1}{N} \sum_{i=1}^N \Big(\frac{1}{{ n}} \sum_{k=1}^{n { - t}} X_i(k) X_i(k+t) \Big)^2 \to_{\rm p} \E \Big( \frac{a^{t}}{1-a^2} \Big)^2.
\end{equation}
By the LLN,
$\frac{1}{N} \sum_{i=1}^N ( \frac{a^t_i}{1-a^2_i} )^2 \to_{\rm p} \E ( \frac{a^t}{1-a^2} )^2$. Therefore by
Minkowski's inequality, for \eqref{X1} we only need to show that
\begin{equation*}
\frac{1}{N} \sum_{i=1}^N \Big( \frac{1}{{ n}} \sum_{k=1}^{n} X_i(k) X_i(k+t) - \frac{a_i^t}{1-a_i^2} \Big)^2  = o_{\rm p}(1).
\end{equation*}
By taking expectations this follows from
\begin{equation}\label{X2}
\E \Big( \frac{1}{n} \sum_{k=1}^{{ n}} X_i(k) X_i(k+t) - \frac{a_i^t}{1-a_i^2} \Big)^2 = \frac{1}{n^2} \E \operatorname{var}
\Big[ \sum_{k=1}^{{ n}}X_i(k) X_i(k+t)  \Big| \Big. a_i \Big] = o(1),
\end{equation}
Using ${\rm cov}[X_i(k) X_i(k+t), X_i(k') X_i(k'+t)|a_i] = {  \frac{a^{2(|k-k'|+t)}}{1-a^4} {\rm cum}_4 + \frac{a^{2|k-k'|} + a^{2\max\{|k-k'|,t\}}}{(1-a^2)^2}}$ and the { same bound as} in
\eqref{var1} we see that the l.h.s. of \eqref{X2} does not exceed $C \E [\frac{1}{(1-a_i)^2} \min \{ 1, \frac{1}{n(1-a_i)}\} ] $
which vanishes as $n \to \infty$ by the { dominated convergence theorem},   due to $\E (1-a)^{-2}<\infty$. \hfill $\Box$
}

\bigskip

\section*{Acknowledgments}

The authors are grateful to an anonymous referee and associate editor for useful comments.
Vytaut{\.e} Pilipauskait{\.e} acknowledges the financial support from the project ``Ambit fields: probabilistic properties and statistical inference'' funded by Villum Fonden.

\end{document}